\documentclass[11pt,a4paper,fleqn]{cas-sc}
\usepackage[natbibapa]{apacite}
\usepackage{gensymb}
\usepackage[section]{placeins}
\usepackage{algpseudocode}
\usepackage{float}
\usepackage{wasysym}

\usepackage{multirow}

\usepackage[linesnumbered,lined,boxed,commentsnumbered,ruled,longend]{algorithm2e}
\usepackage{amsmath,mathtools}

\newcommand{\RN}[1]{%
  \textup{\uppercase\expandafter{\romannumeral#1}}%
}

\begin{document}

\let\WriteBookmarks\relax
\def\floatpagepagefraction{1}
\def\textpagefraction{.001}
\shorttitle{Benders decomposition under climate uncertainty}
\shortauthors{Göke, Schmidt, and Kendziorski}

\title [mode = title]{Stabilized Benders decomposition for energy planning under climate uncertainty}                  

\author[1,2,3]{Leonard Göke}[type=editor,auid=000,bioid=1]
\cormark[1]
 \ead{lgoeke@ethz.ch}
 \author[2]{Felix Schmidt}[auid=000,bioid=1]
\author[1,2]{Mario Kendziorski}[auid=000,bioid=1]

\affiliation[1]{
             organization={Workgroup for Infrastructure Policy (WIP), Technische Universität Berlin},
             addressline={Straße des 17.\,Juni 135}, 
             city={10623 Berlin},
             country={Germany}}
 \affiliation[2]{
             organization={Energy, Transportation, Environment Department, German Institute for Economic Research (DIW Berlin)},
             addressline={Mohrenstraße 58}, 
             city={10117 Berlin},
             country={Germany}}
\affiliation[3]{
             organization={Energy and Process Systems Engineering, Department of Mechanical and Process Engineering, ETH Zurich},
             addressline={Tannenstrasse 3}, 
             city={Zurich 8092},
             country={Switzerland}}

\cortext[cor1]{Corresponding author.}

\begin{abstract}
This paper applies Benders decomposition to two-stage stochastic problems for energy planning under climate uncertainty, a key problem for the design of renewable energy systems. To improve performance, we adapt various refinements for Benders decomposition to the problem's characteristics---a simple continuous master-problem, and few but large sub-problems. The primary focus is stabilization, specifically comparing established bundle methods to a quadratic trust-region approach for continuous problems.

An extensive computational comparison shows that all stabilization methods can significantly reduce computation time. However, the quadratic trust-region and the non-quadratic box-step method are the most robust and straightforward to implement. When parallelized, the introduced algorithm outperforms the vanilla version of Benders decomposition by a factor of 100. In contrast to off-the-shelf solvers, computation time remains constant when the number of scenarios increases.

In conclusion, the algorithm enables robust planning of renewable energy systems with a large number of climatic years. Beyond climate uncertainty, it can make an extensive range of other analyses in energy planning computationally tractable, for instance, endogenous learning and modeling to generate alternatives.
\end{abstract}


\maketitle

\section{Introduction}

The global energy system must undergo a major transformation and replace 86\% of primary energy from fossil fuels as of 2019 with renewable resources to achieve the objectives of the Paris Climate Agreement and limit global warming, \citep{owidenergya}. At the heart of this transformation is renewable electricity from wind and photovoltaic (PV) for two reasons: First, its technical potential exceeds other renewables and even demand projections \citep{Creutzig2017}. Estimates for the global potential of PV alone range from 1,585 to 50,580\,EJ---at least triple the primary consumption in 2019. Second, the use of renewable electricity is not limited to the power sector but can also be deployed for electric heating, mobility, and the creation of synthetic fuels to decarbonize the heat, transport, and industry sectors.

\subsection{System planning under climate uncertainty}

Techno-economic planning models are critical tools for analyzing the transformation towards renewable energy systems. Models decide on the expansion and operation of technologies based on representing the flow and conversion of energy in future systems. Commonly, models are linear optimization problems that minimize the costs to satisfy an exogenous demand considering different boundary conditions, for instance, emission limits. Typical applications include developing long-term scenarios for the energy system, assessing technologies in a system context, and analyzing energy policy \citep{Goeke2021}.

The fluctuating nature of PV and wind generation is challenging for system planning. In conventional energy systems dominated by dispatchable thermal plants, reliable planning only requires a small number of representative time-steps. However, renewable systems require much higher temporal resolution to capture fluctuations of PV and wind. Additionally, planning must consider energy storage as a key option to balance fluctuations \citep{Levin2023}. Both characteristics significantly increase problem size and complexity, quickly rendering the underlying linear optimization problem computationally intractable \citep{Goeke2021a}. As a result, the literature proposes different techniques to reduce model size while representing renewables accurately, for instance, iteratively adjusting the representative time-periods or limiting high resolution to selected parts of the system \citep{Teichgraeber2021,Goeke2020b}.

Nevertheless, all methods for tractable planning of renewable energy systems can only consider one representative year of climate conditions, but research shows that renewable supply and energy demand vary substantially across years \citep{Pfenninger2016,Pfenninger2017}. Based on a sample of 40 historical years, \citet{Bloomfield2016} investigate effects on the British power system and conclude robust planning should consider at least ten of these years. \citet{Ruhnau2022} study the effect of inter-annual variability on storage requirements for a stylized fully renewable German power system. Compared to a representative year, storage requirements double when considering 35 years of climate data. \citet{Ohlendorf2020} specifically analyze the frequency of low-wind power events in Germany threatening system adequacy, again concluding that planning models should consider more climatic years. 

Overall, climate uncertainty and risk-aware planning are essential for the security of supply in renewable energy systems. Existing heuristic approaches to consider climate uncertainty are limited: 
Repeatedly solving the same planning problem for different climate conditions can indicate the extent of variability but not provide a generally reliable solution \citep{Grochowicz2023}. Solving a single problem with representative time-steps from a multi-year sample of climate conditions cannot capture the entire range of climate conditions \citep{Hilbers2019}. Accordingly, recent publications call for new planning methods that deploy stochastic optimization to consider the uncertainty of climate conditions, especially as climate change further increases uncertainty \citep{Craig2022,Plaga2023,DossGollin2023}. 

Considering climate uncertainty in system planning creates a two-stage stochastic problem: The first stage decides on capacity expansion with uncertainty regarding climate conditions. The second stage reveals the climate conditions and decides on the operation of capacities, e.g., the production of a thermal power plant. Stochastic programming uses independent scenarios with distinct probabilities to represent different climate conditions in the second stage.

So far, two-stage energy planning is limited to a few scenarios representing a short representative period of a week at maximum \citep{Backe2022}. This approach limits problem size and keeps the model tractable, but it is unsuitable for planning renewable systems. Here, scenarios must be long and cover at least one year of climate conditions to model seasonal storage. In addition, scenarios must cover various climate conditions to capture fluctuations and enable reliable planning. The resulting two-stage problem with multiple one-year scenarios is a simple but large linear problem that off-the-shelf solvers can hardly solve.

\subsection{Benders decomposition for two-stage stochastic problems}

Benders decomposition (BD), first introduced in \citet{Benders1962}, is a decomposition technique to solve optimization problems. \citet{Slyke1969} first applied a variation of BD, termed the L-shaped method, to solve two-stage stochastic problems and potentially improve their computational tractability. In energy planning under climate uncertainty, BD decomposes the problem into a master-problem (MP) addressing technology expansion and several mutually independent sub-problems (SPs) for each scenario. Afterward, the MP and SPs are solved repeatedly to generate constraints, so-called Benders cuts, that are added to the MP until the algorithm converges \citep{Conejo2006}. 

In existing applications of BD in energy planning models, each SP corresponds to a short time-period, in some models as short as an hour \citep{Skar2014,Lohmann2016,Brandenberg2021,Jacobson2023}. As already discussed, such decomposition is not applicable when modeling renewable systems dependent on energy storage. In this case, each SP must cover at least one year to account for seasonal storage. Combined with high temporal detail, this results in a simple MP and several large SPs. Furthermore, MP and SPs can contain discrete and continuous variables. However, in the MP, at least renewable expansion is continuous, considering the size of a single wind turbine or PV panel is small from a macro perspective.

For many problems, the original BD converges slowly and is not competitive to off-the-shelf solvers for the deterministic equivalent. However, extensive research proposes different enhancements to improve the original BD. In a review, \citet{Rahmaniani2017} group them into the following four most important groups, with three of them being relevant to this paper:
\begin{itemize}
  \item \textbf{Solution procedure}: The majority of enhancements in this area reduce the computation time of the MP, like removing non-binding cuts or not solving to optimality in early iterations. In our case, the more relevant SPs can deploy distributed parallel computing. Furthermore, not solving SPs to optimality can produce inexact but valid cuts \citep{Zakeri2000,deOliviera2020}. Alternatively, only selected SPs are solved in each iteration to reduce the computational load \citep{van_ackooij_incremental_2018}.
  
  \item \textbf{Solution generation}: Again, most enhancements regarding solution generation aim to accelerate the MP, for instance, by solving a linear relaxation instead of a mixed-integer problem. Other adjustments to the MP improve the convergence rate, rendering them sensible even if the MP is small compared to the SPs. These include:

  \begin{itemize}
  \item \textit{Multi-cut reformulations} that use separate cuts for each SP instead of a single aggregated cut, improving convergence but making the MP harder. Multi-cut reformulations generally outperform aggregated cuts for energy planning problems since the MP is easy \citep{Jacobson2023}.
  \item \textit{Problem-specific heuristics} that provide an initial feasible solution for BD to obtain an upper bound on the objective value or exclude infeasible solutions at an early stage. To achieve the latter, additional constraints, known as valid inequalities, which are guaranteed to hold in the optimum, are added to the MP, narrowing down its solution space \citep{Cordeau2006}. For instance, \citet{Wang2013} add the second-stage constraints of one scenario to the MP, improving convergence.
  
  \item \textit{Regularization} restricting the MP to solutions close to a reference point, usually the current best solution, to avoid heavy oscillation \citep{Ruszczynski2003}. These methods are not limited to BD or stochastic programming but have their origins in general convex optimization. For instance, \citet{Linderoth2003} apply the $l_\infty$-norm to limit the maximum difference across all variables between a new and the current best solution of the MP; a method first introduced in a general optimization context in \citet{Marsten1975}. Level and proximal bundle methods originate from nonsmooth convex optimization and add a quadratic term to the objective function for stabilization \citep{Ruszczynski1986,Lemarechal1981}. A series of publications applies explicitly these methods to mitigate the oscillation of BD \citep[e.g.][]{Ackooij2014,Zverovich2012}.
  \end{itemize} 

  \item \textbf{Cut generation}: \citet{Magnanti1981} note that degenerate SPs, e.g., SPs without unique solutions, can generate a set of different cuts, all affecting convergence differently. To select strong cuts that improve convergence, they propose to solve a modified version of the dual SP after the original SP. As a result, the method is most efficient if the MP is difficult and the SPs are small.

\end{itemize}

In conclusion, many established enhancements aim at problems with a difficult MP and small SPs, typically combinatorial problems with discrete complicating variables. \citet{Rahmaniani2017} correspondingly state \textit{"it has often been reported that more than 90\% of the total execution [...] time is spent on solving the MP"}. As a result, applicability to the two-stage planning problem investigated in this paper with a simple MP but large SPs is limited. In contrast to typical applications of BD, it is not combinatorial complexity but the size of the individual SPs that makes our problem computationally challenging.

\subsection{Contribution}

In this paper, we apply BD to solve a two-stage stochastic problem with discrete scenarios used for planning renewable macro-energy systems under climate uncertainty. The existing literature neither addresses this problem specifically nor generally covers BD for continuous linear problems (LPs) that are challenging due to the size of the individual subproblems.

The main contribution of the paper is to adapt and apply existing refinements of BD to the given problem. Overall, our contribution is rather practical than theoretical and has a focus on efficiently solving the problem at hand. The examined refinements can be divided into two groups.
\begin{enumerate}
\item First, common improvements for BD including inexact cuts, valid inequalities, and parallelization. Many other state-of-the-art refinements are not examined since they do not suit the specific structure of the problem.
\item Second, (quadratic) stabilization, or regularization, methods that originate from nonsmooth convex optimization and generally address the instability of cutting plane methods. A particular emphasis is on comparing established methods with a quadratic trust-region approach, which has not been widely used for BD so far.
\end{enumerate}

The structure for the remainder of the paper is as follows: The next section \ref{base} generally introduces the investigated planning problem and a standard implementation of the Benders algorithm. The following section \ref{ref2} focuses on different stabilization methods and corresponding parameter strategies. Afterwards, section \ref{algRef} discusses additional refinements. Section \ref{case} introduces a specific planning problem for an extensive computational comparison of algorithm setups in section \ref{bench}. The final section concludes and suggests further developments.

\section{Benders algorithm} \label{base}

The two-stage energy planning problem decides on capacity expansion in the first stage and operation in the second. Operation is subject to uncertainty represented by different scenarios, for instance, corresponding to different climatic years. The following subsections first introduce the closed formulation of the problem, followed by the basic BD implementation used to solve it. The notation uses lowercase letters for parameters, uppercase letters for variables, and Greek letters for anything specific to BD.

\subsection{Problem formulation} \label{form}

Eqs. \ref{eq:1a} to \ref{eq:1h} provide the closed formulation of the two-stage stochastic energy planning problem. For a set of technologies $i \in I$, the problem decides on expanding capacities $Exp_{y,i}$ over a set of chronological but not necessarily consecutive years $y \in Y$. All variables of the problem are continuous and non-negative.

According to \ref{eq:1b}, the capacity $Capa_{y,i}$ available in each year depends on expansion in the current and previous years provided by $Y_y^{exp}$. In Eqs. \ref{eq:1c} to \ref{eq:1f} these capacities constrain the operation of technologies at each time-step $t \in T$, in each scenario $s \in S$, and within each year. For generation technologies $I^{ge}$, the capacity constraint limiting the generation $Gen_{y,s,i,t}$ includes a scenario-dependent capacity factor $cf_{s,i,t}$ reflecting the available share of capacity. For storage technologies, $I^{st}$, two different capacities constrain three different operational variables: Charged energy $St_{y,s,i,t}^{in}$ and discharged energy $St_{y,s,i,t}^{out}$ are both constrained by the power capacity $Capa_{y,i}^{st}$; the current storage level $St_{y,s,i,t}^{size}$ is constrained by the energy capacity $Capa_{y,i}^{size}$. Net supply from generation and storage technologies must match the exogenous demand $dem_{y,s,t}$ in each year, time-step, and scenario, as expressed in the energy balance in Eq. \ref{eq:1g}. The balance includes a loss-of-load variable $Lss_{y,s,t}$ to avoid an infeasible problem. The storage balance in Eq. \ref{eq:1b} computes the storage level at time-step $t$ based on the storage level in the previous time-step $t-1$ plus charged and minus discharged energy. The storage constraint is circular for each year, meaning the last time-step in a year is previous to the first.

The objective function of the problem in Eq. \ref{eq:1a} minimizes total system costs comprised of expansion costs $U$ and the sum of operational costs $V_{y,s}$ across all years $y$ and scenarios $s$ weighted according to the scenario probabilities $p_s$. Expansion costs defined in Eq. \ref{eq:1i} depend on the expansion variable and the specific expansion costs $c^{inv}_{y,i}$ of each technology. Operational costs defined in Eq. \ref{eq:1j} depend on the costs associated with loss-of-load $c^{lss}$ and variable costs of generation technologies $c^{var}_{y,s,i}$, for instance, fuel costs.
{\small
\begin{subequations}
\begin{alignat}{4}
\min  & \; \;  U  && + && \sum_{y \in Y, s \in S} p_{s} \cdot V_{y,s} && \label{eq:1a} \\
\text{s.t.} \, \, &  U  && \, = \, &&  \sum_{y \in Y, i \in I} c^{inv}_{y,i} \cdot Exp_{y,i}  && \label{eq:1i}  \\
& V_{y,s}  && \, = \, && \sum_{t \in T} Lss_{y,s,t} \cdot c^{lss} + \sum_{i \in I} Gen_{y,s,i,t} \cdot c^{var}_{y,s,i} && \, \forall y \in Y, s \in S \label{eq:1j} \\
& Capa_{y,i} && \; = \; && \sum_{y' \in Y_{y}^{exp}} Exp_{y',i} && \; \forall y \in Y, i \in I \label{eq:1b} \\ 
& Gen_{y,s,i,t} && \, \leq \, && cf_{s,i,t} \cdot Capa_{y,i}^{ge} &&   \, \forall y \in Y, s \in S, i \in I^{ge}, t \in T \label{eq:1c}  \\
& St_{y,s,i,t}^{out} && \, \leq \, && Capa_{y,i}^{st} &&   \, \forall y \in Y, s \in S, i \in I^{st}, t \in T  \label{eq:1d} \\
&  St_{y,s,i,t}^{in} && \, \leq \, && Capa_{y,i}^{st} &&   \, \forall y \in Y, s \in S, i \in I^{st}, t \in T \label{eq:1e}  \\
& St^{size}_{y,s,i,t} && \, \leq  \, && Capa_{y,i}^{size}  &&   \, \forall y \in Y, s \in S, i \in I^{st}, t \in T  \label{eq:1f} \\
& dem_{y,s,t} &&  \, = \, && Lss_{y,s,t} + \sum_{i \in I^{ge}} Gen_{y,s,i,t} - \sum_{i \in I^{st}} St_{y,s,i,t}^{in} + St_{y,s,i,t}^{out} && \,\forall y \in Y, s \in S, t \in T \label{eq:1g} \\ 
& St^{size}_{y,s,i,t} &&  \, = \, && St^{size}_{y,s,i,t-1} + St_{y,s,i,t}^{in} - St_{y,s,i,t}^{out} && \, \forall y \in Y, s \in S, i \in I^{st}, t \in T \label{eq:1h}
\end{alignat}
\end{subequations}
}
The formulation in Eqs. \ref{eq:1a} to \ref{eq:1h} introduces the general problem structure and all details pivotal for BD but is still stylized. The model introduced in section \ref{case} and used for benchmarking includes additional elements. First, the model includes multiple regions of expansion and operation. The model can build transmission capacity between regions in the first stage to exchange energy in the second. Next, distinct energy carriers, like electricity and hydrogen, are converted into one another by respective technologies, like hydrogen turbines or electrolyzers. Finally, there are storage losses and additional constraints, imposing upper limits on capacities or restricting operation, for instance, a yearly emission limit. For a comprehensive formulation of the planning problem, see \citet{Goeke2020b}.

\subsection{Benders decomposition} \label{bdStd}

\begin{figure}[!htp]
	\centering
		\includegraphics[scale=0.50]{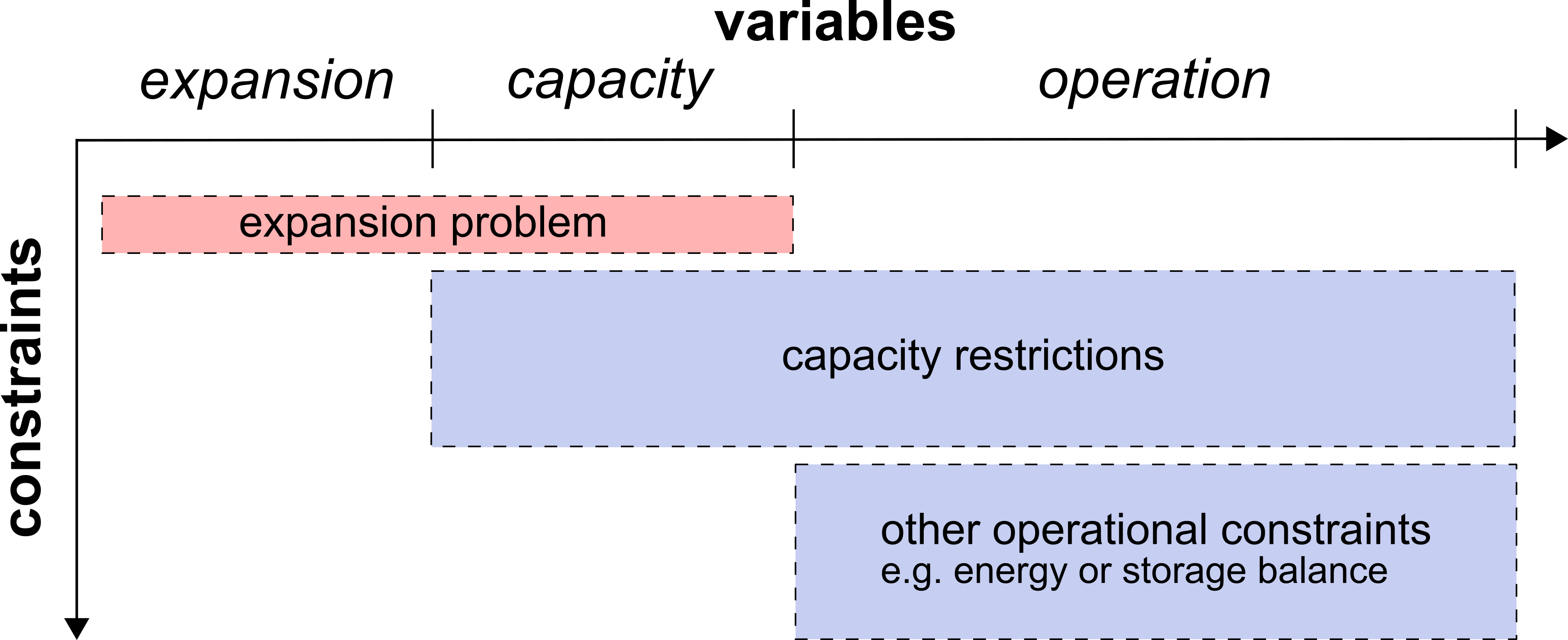}
	\caption{Structure of closed two-stage problem}
	\label{fig:baseProb}
\end{figure}

Fig. \ref{fig:baseProb} is an illustrative depiction of the planning problem based on the block structure of the constraint matrix. The expansion problem in the first stage includes expansion and capacity variables, the constraint in \ref{eq:1b} connecting them, and the cost definition in Eq. \ref{eq:1i}. The second stage includes the remaining constraints, namely the energy and storage balance, the definition of variable costs, and the capacity restrictions. The latter links capacity with operational variables and connects the first stage and second stage of the problem. 

Even for a single year and scenario, the number of constraints and variables of the expansion problem in the first stage is small compared to the operation in the second stage. Increasing the number of years extends both stages, but increasing the number of stochastic scenarios extends the operational problem only. Since the expansion problem is small, the total size almost linearly depends on the number of stochastic scenarios, and the problem quickly becomes intractable if their number increases.

\begin{figure}[!htp]
	\centering
		\includegraphics[scale=0.50]{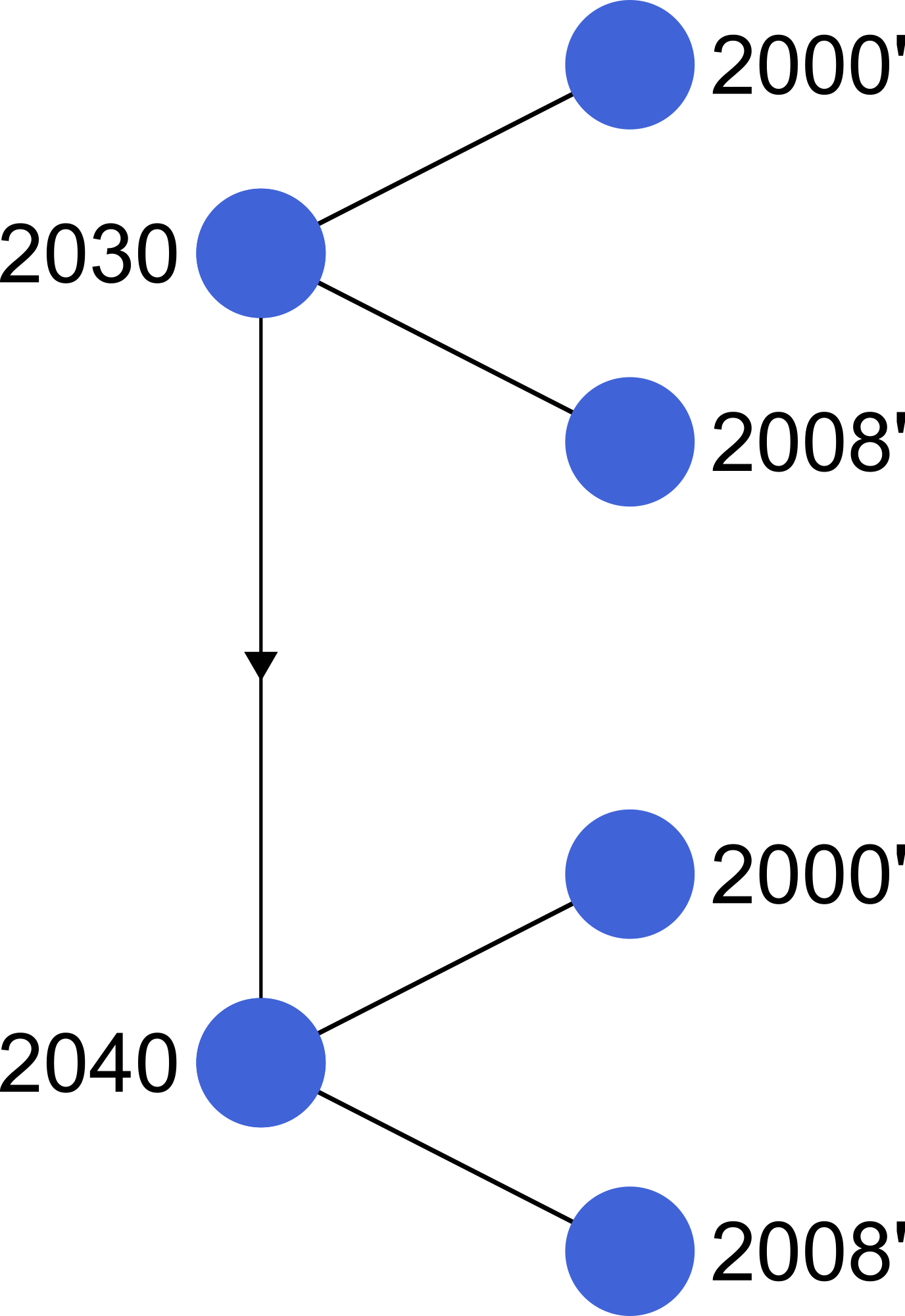}
	\caption{Graph of years and scenarios for exemplary problem}
	\label{fig:tree}
\end{figure}

For BD, the first stage with expansion corresponds to the MP, while several independent SPs cover operation in the second stage. Consider the temporal problem structure outlined in Fig. \ref{fig:tree} for demonstration. The first stage models two steps of capacity expansion with a 10-year interval. Multi-year intervals are a common modeling strategy for representing transformation pathways, aiming to reduce computational complexity. The second stage includes two operational scenarios for each expansion step. These scenarios represent different patterns of renewable supply and energy demand. In this example, these patterns use historical data for the climatic years $2000'$ and $2008'$. These two years are selected at random here for illustrative purposes. 
For details on the methods for selecting a sample of climatic years, see section \ref{clus} in the Appendix. In addition, it is conceivable to have different operational scenarios for each expansion step, for instance, to account for changes in patterns caused by climate change.

According to the problem formulation in the previous section, operation is modeled for each year and scenario, resulting in four distinct operational problems, as shown in Fig. \ref{fig:tree}. Each operational problem can be solved independently for a given set of capacities if no complicating constraint links the variables from different SPs. An example of such a constraint is an emission limit not enforced for each SP separately but jointly across different years and scenarios, like a carbon budget.

\begin{figure}[!htp]
	\centering
		\includegraphics[scale=0.50]{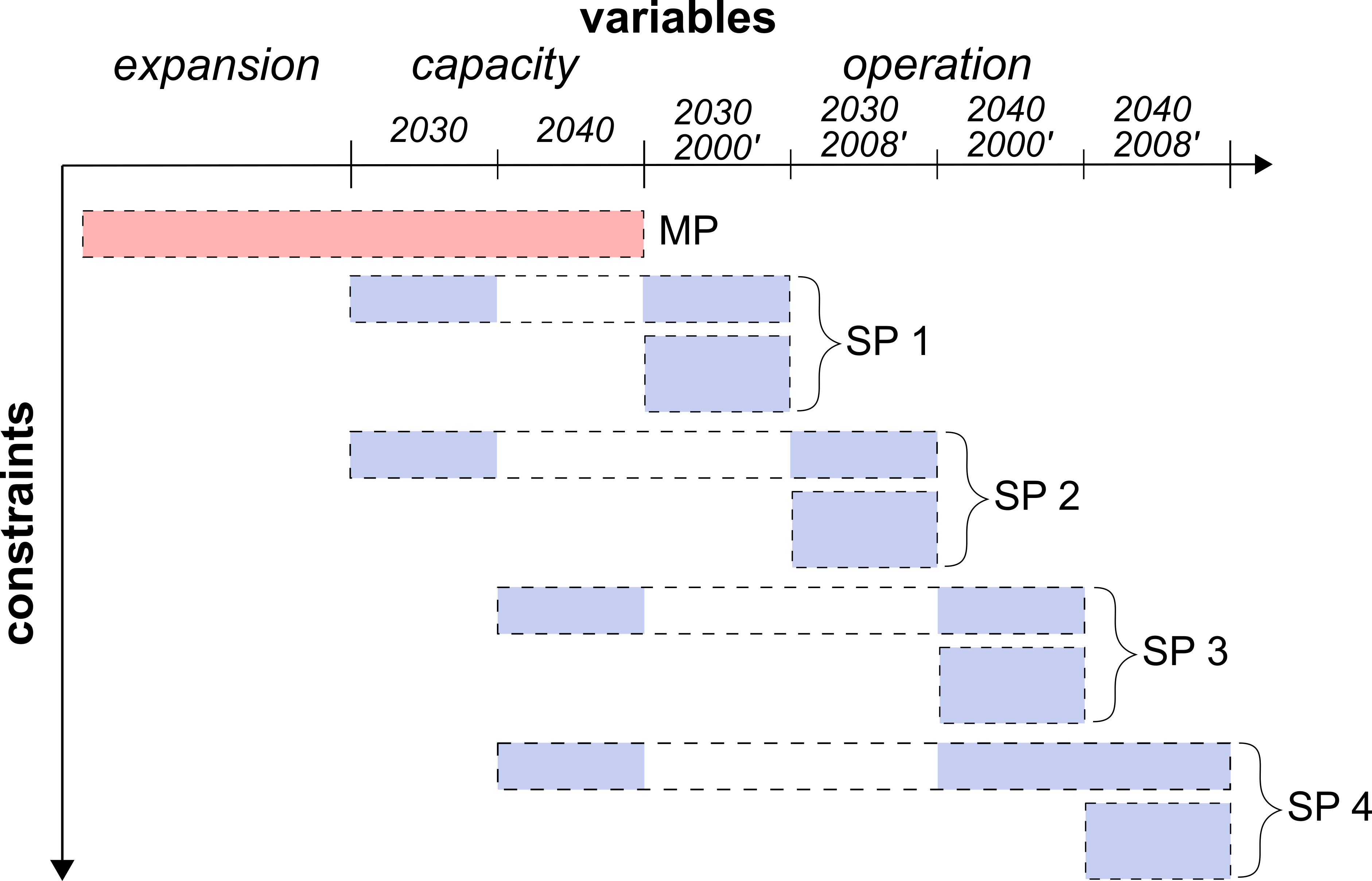}
	\caption{Structure of decomposed problem}
	\label{fig:matrix}
\end{figure}

Fig. \ref{fig:matrix} provides the structure of the decomposed problem, analogously to Fig. \ref{fig:baseProb}. It shows how the approach splits the operational problem into four independent SPs, each consisting of different operational variables but sharing capacities with the MP and other SPs for the same year. Splitting the second stage operation into several SPs is beneficial because it enables distributed computing and prevents complexity in the second stage to scale with the number of scenarios. In addition, it increases the number of cuts added to the MP in each iteration, improving convergence, as we will elaborate when introducing the BD algorithm next. 

Alg. \ref{alg:stdBD} presents the standard version of the BD algorithm applied in this paper. $k$ is an iteration counter; $\zeta_{low}$ and $\zeta_{up}$ are the objective function's current lower and upper bound, respectively. We use a multi-cut version and, at each iteration, add separate constraints for each sub-problem to the MP instead of a single aggregated cut. The approach improves convergence but increases the complexity of the MP---a favorable trade-off for linear planning probelms since the MP is simple \citep{Jacobson2023}.

\begin{algorithm} \label{alg:stdBD}
\small
\caption{Multi-cut benders algorithm}
\DontPrintSemicolon
\SetKwBlock{DoParallel}{do in parallel}{end}
\textbf{\textit{Step 0 (Initialization)}}\;
chose convergence tolerance $\epsilon$ and deletion threshold $\eta$ \;
set $k = 1$, $\zeta_{up} \leftarrow \infty$ and $\zeta_{low} \leftarrow 0 $ \;
\textbf{\textit{Step 1 (Convergence test)}}:\;
\textbf{if} $1 - \frac{\zeta_{low}}{\zeta_{up}} < \epsilon$ \textbf{then} \textit{stop} \textbf{end if}\;
\textbf{\textit{Step 2 (Obtain new iterate)}}\;
solve $\textbf{MP}$ 
get $Capa^{(k)}$ and $U^{(k)}$ \;
set $\zeta_{low} \gets obj(\mathbf{MP})$ \;
\DoParallel{
\For{$y \in Y, \, s \in S$}{
\textbf{\textit{Step 3.1 (Solve sub-problems)}}\;
solve $\textbf{SP}_{y,s}$ with $Capa^{(k)}$ \;
get $V_{y,s}^{(k)}$ and $\lambda_{y,s}^{(k)}$ \;
\textbf{\textit{Step 3.2 (Add cuts to MP)}}\;
add $\Omega_{k,y,s}$ to \textbf{MP} \;
set $\delta_{k,y,s} \gets k$ \;
}
}
\textbf{\textit{Step 4 (Bundle management)}}\;

\For{$l \in \{1,\dots,k-1\}, y \in Y, \, s \in S$}{
\uIf{$\Omega_{l,y,s}$ is binding}{
set $\delta_{l,y,s} \gets k$ \;
}
\uElseIf{$k - \delta_{l,y,s} > \eta$}{
delete $\Omega_{l,y,s}$ from \textbf{MP} \;
}
}
\textbf{\textit{Step 5 (Descent test)}}\;
\If{$\zeta_{up} >  U^{(k)} + \sum_{y \in Y, s \in S} p_{s} \cdot V_{y,s}^{(k)}$}{
set $\zeta_{up} \gets U^{(k)} + \sum_{y \in Y, s \in S} p_{s} \cdot V_{y,s}^{(k)}$ \;
}
set $k \gets k + 1$ \;
\textbf{return to} \textbf{\textit{Step 1}}
\end{algorithm}

After initializing all variables, the algorithm starts iterating with a convergence check (Step 1), which tests whether the optimality gap is within a pre-defined tolerance. Step 2 solves the original MP, as described by Eqs. \ref{eq:2a} to \ref{eq:2c}. Since the MP only includes the constraints of the expansion problem, its objective function approximates the actual costs of each SP $V_{y,s}$ using the estimator $\alpha_{y,s}$. After solving, the algorithm obtains the resulting capacity vector $Capa^{(k)}$ and expansion costs $U^{(k)}$ for the current iteration and sets the lower bound $\zeta_{low}$ to the current objective of the MP. In the first iteration, $\alpha_{y,s}$ is still unconstrained, and the trivial optimum for the MP is zero, meaning no capacity expansion at all.

\begin{subequations}
\begin{alignat}{4}
\min  &  \; \;  U   && + &&  \sum_{y \in Y, s \in S} p_{s} \cdot \alpha_{y,s} && \label{eq:2a}  \\
\text{s.t.} \, \, & \; \; U  && \, = \, &&  \, \, \, \, \, \sum_{y \in Y} c^{inv}_{y,i} \cdot Exp_{y,i}  &&  \label{eq:2b} \\ 
&  Capa_{y,i} && \; = \; &&  \, \sum_{y' \in Y_{y}^{exp}} Exp_{y',i} && \; \forall y \in Y, i \in I \label{eq:2c} 
\end{alignat}
\end{subequations}

In Step 3.1, the algorithm solves each SP with capacities fixed to the previously computed $Capa^{(k)}$, as described by Eqs. \ref{eq:3a} to \ref{eq:3i}. We will discuss the parallelization mentioned in the algorithm in section \ref{parallel}. Afterward, the algorithm gets the actual costs of each SP $V_{y,s}^{(k)}$ and the dual values $\lambda_{y,s,i}^{(k)}$ of the constraint fixing capacities in Eq. \ref{eq:3i}. In our case, SPs are always feasible since the loss-of-load $Lss_{y,s,t}$ serves as a slack variable for unmet demand.

{\small
\begin{subequations}
\begin{alignat}{4}
\min  & \; V_{y,s}^{(k)} && && &&  \label{eq:3a} \\
\text{s.t.} \, \, & V_{y,s}^{(k)}  && \, = \, && \sum_{t\in T} Lss_{y,s,t} \cdot c^{lss} + \sum_{i \in I, t \in T} Gen_{y,s,i,t} \cdot c^{var}_{y,s,i} && \, \forall y \in Y, s \in S \label{eq:3b}\\
& Gen_{y,s,i,t} && \, \leq \, && cf_{s,i,t} \cdot Capa_{y,i}^{ge} &&   \, \forall y \in Y, s \in S, i \in I^{ge}, t \in T  \label{eq:3c}\\
& St_{y,s,i,t}^{out} && \, \leq \, && Capa_{y,i}^{st} &&   \, \forall y \in Y, s \in S, i \in I^{st}, t \in T  \label{eq:3d}\\
&  St_{y,s,i,t}^{in} && \, \leq \, && Capa_{y,i}^{st} &&   \, \forall y \in Y, s \in S, i \in I^{st}, t \in T  \label{eq:3e}\\
& St^{size}_{y,s,i,t} && \, \leq  \, && Capa_{y,i}^{size}  &&   \, \forall y \in Y, s \in S, i \in I^{st}, t \in T  \label{eq:3f}\\
& dem_{y,s,t} &&  \, = \, && Lss_{y,s,t} + \sum_{i \in I^{ge}} Gen_{y,s,i,t} - \sum_{i \in I^{st}} St_{y,s,i,t}^{in} + St_{y,s,i,t}^{out} && \,\forall y \in Y, s \in S, t \in T \label{eq:3g}\\ 
& St^{size}_{y,s,i,t} &&  \, = \, && St^{size}_{y,s,i,t-1} + St_{y,s,i,t}^{in} - St_{y,s,i,t}^{out} && \, \forall y \in Y, s \in S, i \in I^{st}, t \in T \label{eq:3h}\\
& Capa_{y,i} && = && Capa^{(k)}_{y,i} \, : \, \lambda_{y,s,i}^{(k)} && \forall i \in I \label{eq:3i}
\end{alignat}
\end{subequations}
}

Step 3.2 computes the Benders cuts $\Omega_{k,y,s}$ defined by Eq. \ref{eq:cut} and adds them to the MP. Each cut improves the estimator $\alpha_{y,s}$ based on the exact result of the SP $V_{y,s}^{(k)}$ and its gradient $\lambda_{y,s,i}^{(k)}$ at the point $Capa_{y,i}^{(k)}$. In other words, by adding hyperplanes restricting the estimator, BD performs a piece-wise linear approximation of the unknown but convex function that assigns capacities to SP-costs \citep{Conejo2006}. In the literature, cuts of this form are termed optimality cuts instead of feasibility cuts that the algorithm adds when the SP is infeasible, which cannot occur in our case. Section \ref{algRef} discusses the bundle management in Step 4 already included in the algorithm here.

\begin{equation}
\Omega_{k,y,s}: \alpha_{y,s} \geq V^{(k)}_{y,s} + \sum_{i \in I} \lambda_{y,s,i}^{(k)} \cdot (Capa_{y,i} - Capa^{(k)}_{y,i}) \label{eq:cut}
\end{equation}

Once all SPs are solved, summing expansion costs $U$ and accurate operational costs $V_{y,s}^{(k)}$ in Step 5 gives the total costs at $Capa^{(k)}$. If these undercut the current $\zeta_{up}$, $Capa^{(k)}$ is a new best solution, and we adjust $\zeta_{up}$ accordingly. Finally, the iteration counter increases by one, and the next iteration start if the convergence test in Step 1 fails.

\section{Stabilization} \label{ref2}

Stabilization is the pivotal refinement of the developed BD algorithm. The stabilization methods discussed below are not restricted to stochastic programming but derive from a much broader field of literature concerned with cutting plane methods. Cutting plane methods iteratively minimize convex functions, which may be non-differentiable, by approximating the function through outer linearizations \citep{Kelley1960}. Consequently, BD is a cutting plane method for two-stage problems. 

In the next section, we describe how the BD algorithm introduced in the previous section can incorporate stabilization. Afterwards, we introduce three groups of methods for stabilization and strategies to parameterize them.

\subsection{Stabilized Benders algorithm} \label{stabAlg}

We begin by representing the BD algorithm above in a more general form to connect it to the literature on stabilization in non-smooth convex optimization. 

Let $\phi_{y,s}(Capa)$ be the optimal objective value of the SP for investment year $y$ and scenario $s$ as a function of complicating variable choice $Capa$. This function and its weighted average are polyhedral (and thus convex and non-smooth), which allows constructing a cutting plane model \citep[Prop. 2.2 \& 2.3]{Shapiro2009}. 

The cutting plane approximation of $\phi_{y,s}(Capa)$ at iteration $k$ corresponds to the minimum of the estimator $\alpha_{y,s}$ restricted by the cuts $\Omega_{k',y,s}$, as shown in Eq. \ref{eq:cutPl}. Each cut, representing a supporting hyperplane of $\phi_{y,s}$, is defined by Eq. \ref{eq:cut}. This formulation is equivalent to the piece-wise maximum of all linearizations across the domain of $Capa$, as the convex optimization literature usually defines the cutting plane model.
\begin{align}
\hat{\phi}_{y,s}^{(k)}(Capa) = \displaystyle{\min_{k > k'}\{\alpha_{y,s} \, | \, \Omega_{k',y,s}\}} \label{eq:cutPl}
\end{align}
The sum of approximations over all SPs $\hat{\phi}^{(k)}(Capa)$ can accordingly be written as Eq. \ref{eq:cutPl2}.
\begin{align}
    \hat{\phi}^{(k)}(Capa) = \sum_{y \in Y, s \in S} p_s \hat{\phi}_{y,s}^{(k)}(Capa) \label{eq:cutPl2}
\end{align}

By sending a candidate solution $Capa^{(k)}$ to an \textit{oracle}, which corresponds to solving the SPs in our case, we obtain $\sum_{y,s} p_s \phi_{y,s}(Capa^{(k)})$ and subgradients $\lambda_{y,s,i}^{(k)}$, which are in the subdifferentials of the respective $\phi_{y,s}(Capa^{(k)})$. Subgradients generalize the concept of gradients to non-differentiable functions. $g$ is a subgradient of $f(x)$ if $f(y)\geq f(x) + \langle g,y-x \rangle \ \forall y$. The subdifferential of $f(x)$ is the set of all subgradients at $x$. The subdifferential is a singleton if the function is differentiable at $x$. On this basis, we improve the cutting plane model in Eq. \ref{eq:cutPl2} by adding new cuts and then solve the MP, $\min U(Capa) + \hat{\phi}^{(k)}(Capa)$, to obtain a new candidate solution $Capa^{(k+1)}$.

Although the proof of convergence for $\min \hat{\phi}(Capa) \to \min \phi(Capa)$ is simple \citep{Kelley1960}, cutting plane algorithms suffer from two computational issues. First, the algorithm is inherently unstable since the capacities $Capa$ in two consecutive iterates can be extreme points of the feasible space arbitrarily far apart. In energy planning, this problem is particularly pronounced because the feasible region is hardly restricted and completely continuous. Even enforcing the optimum as a starting value may not improve the solution time since the algorithm is initially unable to evaluate the quality of the solution. Second, the size of the model approximating $\phi(Capa)$, in the case of BD the number of cuts, can become prohibitively large as the number of iterations $k$ increases and slow the algorithm \citep{Frangioni2020}. We emphasize that these issues are mutually reinforcing: Oscillation increases the number of iterations, which in turn slows down the convergence of the approximating model $\hat{\phi}^{(k)}(Capa)$. 

The term "bundle methods" encompasses methods to address these two issues by managing the "bundle" of information accumulated throughout the iterations \citep{Frangioni2002}. In the following, our exclusive focus will be on stabilization methods that center each iteration around a reference point to prevent oscillation \citep{Bonnans2006}. Since these methods improve convergence at the expense of adding complexity to the MP, they are well-suited to the examined problem with a simple MP but costly oracles, viz., large individual SPs. Removing non-binding cuts after a certain number of iterations addresses the second issue, as further discussed in section \ref{algRef}. Additionally, stabilization reduces the overall number of iterations, thereby mitigating problems associated with excessive growth of the bundle as well.

All stabilzation methods in the literature on (non-smooth) convex optimization build on the cutting plane method and aim to mitigate oscillation by keeping the next iterate $Capa^{(k)}$ relatively close to a \textit{stability center}, a reference point in the (feasible) space of complicating variables. Commonly, this point is part of the sequence of previous iterates $\{Capa^{(j)}\}_{j=1,\dots,k-1}$. In theory, various measures can quantify the distance between the stability center and the next iterate, but the literature typically relies on the Euclidean norm. For the proximal bundle method, discussed below, the literature constructs a "poor man's Hessian" approximates valuable second-order information otherwise absent in the polyhedral setting to justify the Euclidean distance \citep{Frangioni2020}.

In the problem described in Section \ref{form}, the capacity vector $Capa$ is implicitly defined by the selected expansion path (cf. Equation \ref{eq:1b}). In practical applications, many elements of the capacity vector are fixed constants. Therefore, defining the stability components below in terms of the lower-dimensional vector of expansion decisions $Exp$ is computationally advantageous. 

An update of the stability center to the current iterate is called a \textit{serious step}, and an iteration that keeps the current stability center is called a \textit{null step}. Nevertheless, a null step is still beneficial since it improves the cutting plane model locally. The literature offers various methods to distinguish between null and serious steps. For example, a \textit{serious step} can be a step that improves the objective by a specified value compared to the current stability center. In an alternative definition, the ratio between the actual improvement determined by $\phi(Capa)$ and the improvement estimated by $\hat{\phi}(Capa)$ must exceed a certain threshold. In our case, stabilization achieves the best results when the algorithm makes a \textit{serious step} in each iteration that improves the current best solution, even marginally.

Algorithm \ref{alg:bundle} represents our general implementation of stabilization with BD. The following sections will extend this general form with specific stabilization methods and corresponding parameter strategies.

\begin{algorithm}[!htp] \label{alg:bundle}
\small
\caption{General parallel stabilized Benders algorithm}
\DontPrintSemicolon
\SetKwBlock{DoParallel}{do in parallel}{end}
\textbf{\textit{Step 0.1 (Parameter Initialization)}} \\
Set convergence tolerance $\epsilon>0$, deletion threshold $\eta$, $\zeta_{low} \leftarrow 0$, $\zeta_{up} \leftarrow \infty$, set dynamic parameter for chosen bundle method $\in\{\tau_1,\gamma_1,\ell_1\}$, define stability center $Exp^{up} \leftarrow \mathbf{0}$, init. serious step counter $q=0$, set $k=1$\;
\textbf{\textit{Step 0.2 (Stability center initialization)}}\;
\If{\textit{initialize stability center true}}{
    solve \textbf{CP} for $S':={s\in S | p_s = max(p_s)}$; obtain $Exp^{(0)}$\;
    \DoParallel{$\mathbf{\forall y,s}:$ solve \textbf{SP}$_{y,s}$ with $Capa^{(0)}$; obtain $V_{y,s}^{(0)}, \lambda_{y,s,i}^{(0)}$, compute $\Omega_{0,y,s}$ add to stabilized \textbf{MP} and \textbf{MP}; $\delta_{0,y,s} = 0$\;
    }
    $Exp^{up} \leftarrow Exp^{(0)}$ 
}
\textbf{\textit{Step 1 (Convergence test)}}\;
\textbf{if} $1-\frac{\zeta_{low}}{\zeta_{up}}<\epsilon$ \textbf{then} \textit{stop} \textbf{end if}\;
\textbf{\textit{Step 2 (Obtain new iterate)}}\;
Solve stabilized \textbf{MP} and obtain $Capa^{(k)}$, $Exp^{(k)}$ and $U^{(k)}$\;
\DoParallel{
\textbf{\textit{Step 3.1 (Oracle call: solve subproblems)}}\;
$\mathbf{\forall y,s}$: solve \textbf{SP}$_{y,s}$ with $Capa^{(k)}$; obtain $V_{y,s}^{(k)}, \lambda_{y,s,i}^{(k)}$, compute $\Omega_{0,y,s}$ add to stabilized \textbf{MP} and \textbf{MP}; $\delta_{k,y,s} = k$\;
\textbf{\textit{Step 3.2 (Obtain lower bound)}}\;
solve \textbf{MP}, set $\zeta_{low} \leftarrow obj(\textbf{MP})$ and obtain $U_{unstab}^{(k)}$
}
\textbf{\textit{Step 4 (Bundle management)}}\;
\For{$ y,s,l\in\{0,\dots,k-1\}$}{
\uIf{$\Omega_{l,y,s}$ is \textit{binding}}{$\delta_{l,y,s} = k$}
\uElseIf{$k-\delta_{l,y,s} > \eta$}{delete $\Omega_{l,y,s}$ from \textbf{MP} and stabilized \textbf{MP} \;
}
}
\textbf{\textit{Step 5 (Descent test)}}\;
\eIf{$U^{(k)} + \sum_{y,s} p_sV_{y,s}^{(k)} < \zeta_{up}$}
{
$\zeta_{up} \leftarrow U^{(k)} + \sum_{y,s} p_sV_{y,s}^{(k)}$ \;
$Exp^{up}\leftarrow Exp^{(k)}$ \;
$q = q+1$ \;
}{
$q = 0$ \;
}
\textbf{\textit{Step 6 (Dynamic parameter adjustment)}}:\;
Follow steps of adjustment algorithms below \;
$k \leftarrow k + 1$
\textbf{Return to \textit{Step 1}}
\end{algorithm}

Step 0.2 of the algorithm produces an initial feasible solution. In this study, we solve a deterministic instance of the model with the most probable scenario. Section \ref{bench} goes more in-depth on testing different initialization strategies. Step 2 obtains the next iterate by solving the stabilized MP, but the algorithm still solves the original MP in each iteration to obtain a lower bound. The variable $q$ counts the number of serious steps and is an input to the Parameter strategies discussed below. Again, to be concise, we show a parallelized version of the algorithm here but discuss the parallelization in section \ref{parallel}.

We did not include a minimum descent target since tests showed that the algorithm does not profit from this method. However, replacing the condition of Step 5 with the equation below can easily implement this feature:
    \begin{align*}
        U^{(k)} + \sum_{y,s} p_sV_{y,s}^{(k)} < \zeta_{up} - \omega_k
    \end{align*}
    where $\omega_k$ gives the descent target in iteration $k$. We refer to \cite{Frangioni2020} for an introductory discussion of a commonly used Armijo-type rule.

\subsection{Bundle methods} \label{bmeth}

There are three relevant bundle methods for the stabililization of cutting plane algorithms and each methods features at least one dynamic parameter that governs the restrictiveness of the stabilization. While all presented methods are theoretically equivalent in the sense that there exists a combination of respective parameter choices that results in the exact next iterate from each bundle method (see \cite{Bonnans2006}, Thm. 10.7 for a proof), their practical performance hinges critically on the definition of a parameter strategy, also referred to as \textit{t}-strategy \citep{Frangioni2002}). Finding good \textit{t}-strategies is not trivial, and selecting a poor one may even produce inferior results compared to an unstabilized cutting plane algorithm. For instance, an overstabilized MP, will not benefit from the global nature of the cutting plane model and may take more iterations than an unstabilized version (cf. \cite{Frangioni2020}, Fig. 3.2). Strategies suggested in the literature are usually comprised of a set of rules that determines the parameter values for the next iteration dependent on a set of variables describing the algorithm's progress. For example, such a rule could prescribe that the dynamic parameter doubles after \textit{l} null steps and halves after \textit{m} serious steps. In this simple strategy \textit{l} and \textit{m} are \textit{hyperparameters} subject to calibration.

In this section, we describe the bundle methods and corresponding parameter strategies. Parameter strategies depend on the method and replace Step 6 in the previously outlined algorithm. Finding good parameter strategies is complex and can be highly problem-specific (\cite{Frangioni2020}). We highlight the hyperparameters used in each strategy before testing specific values in section \ref{bench}. 

\subsubsection{Proximal-bundle methods}

The first bundle method is called \textit{proximal bundle method} (PBM)\citep{lemarechal1977}. This variant adds a penalty term to the MP's objective function to prevent moving away from the stability center. Two-stage stochastic problems widely use this method, for example in \cite{oliviera2011} or \cite{Zverovich2012}. The stabilized MP is a quadratic program and is given by:
\begin{align} \label{eq:prox}
    \min_{\{Capa,Exp\}} U(Capa) + \hat{\phi}^{(k)}(Capa) + \frac{1}{2\tau_k}\|Exp-Exp^{up}\|_2^2 \quad s.t. \ \ Capa_{y,i} = \sum_{y'\in Y_y} Exp_{y',i} \ \forall y,i
\end{align}

However, PBM can go beyond the $l_2$norm. \cite{pessoa_automation_2018} use a generalized linear PBM approach using a composite penalty term consisting of three terms. This approach does not penalize small deviations determined by the $l_\infty$-norm, but penalizes deviations outside this confidence interval based on the $l_1$-norm. As the penalty parameter goes to infinity, this approach turns into the Boxstep method we will introduce in section \ref{qtr}. To keep the complexity of methods comprehensible, we limited our computational study to the standard PBM method.

\paragraph{Parameter strategies}

The dynamic parameter $\tau_k$ in Eq. \ref{eq:prox} in PBM scales the penalty in the objective function with larger values decreasing it; and smaller values increasing it. 

Following \cite{Oliveira2016}, we use two ruled-based parameter strategies first proposed by \cite{kiwiel_proximity_1990}. We implement the strategy by inserting the steps of Algorithm \ref{alg:pbm} into Step 6 of Algorithm \ref{alg:bundle}. Both strategies decrease the weight in response to a serious step and increase it in response to a null step (Step 6.3). However, the methods compute the auxiliary parameter $\tau_k^{aux}$ (Step 6.2) differently. For both methods, the auxiliary parameter mimics the Hessian used in a Newton-type algorithm on a Moreau envelope of the true function $\phi$. Of course, this function is unknown and we can only approximate the the optimal descent. For details, we refer to \cite{Frangioni2020} and \cite{HiriartUrruty1993}, Chapter XV. 

In PBM-1, the auxiliary parameter uses the ratio between the achieved descent and the descent predicted by the model (\cite{kiwiel_proximity_1990}). PBM-2 derives directly from a quasi-Newton formula and requires some safeguards to ensure that the numerator of the fraction term remains non-negative \citep{lemarechal_variable_1997}. Both methods require the following hyperparameters: the starting value $\tau_1$, the minimum parameter value $\tau_{min}$ and the scaling factor $a$.

\begin{algorithm} \label{alg:pbm}
    \small
    \caption{Dynamic adjustment of PBM parameter $\tau_k$}
    \DontPrintSemicolon
    \textbf{\textit{Step 6.1 (Initialization)}}:\;
    initialize hyperparameters $a,\tau_{min},\tau_1$ \;
    \textbf{\textit{Step 6.2 (Define auxilary parameter)}}:\;
    $\tau_k^{aux} \leftarrow 2\tau_k \left( 1+\frac{U(Capa^{up}) - U(Capa^{(k+1)})+\phi(Capa^{up}) - \phi(Capa^{(k+1)})}{U(Capa^{up}) - U(Capa^{(k+1)})+\phi(Capa^{up}) - \hat{\phi}^{(k)}(Capa^{(k+1)})} \right)$ \textbf{if} method is \textit{PBM-1} \;
    $\tau_k^{aux} \leftarrow \tau_k\left(1+\frac{\sum_{y,s}p_s\left[\sum_{i}(\lambda_{y,s,i}^{(k+1)}-\lambda_{y,s,i}^{(k)})\cdot (Capa^{k+1}_{i}-Capa^{k}_{i})\right]}{\sum_{y,s}p_s\left[\sum_{i}(\lambda_{y,s,i}^{(k+1)}-\lambda_{y,s,i}^{(k)})^2\right]}\right)$ \textbf{if} method is \textit{PBM-2} \;
    \textbf{\textit{Step 6.3 (Update dynamic parameter based on step type)}} \;
    \eIf{If $q=0$}
    {$\tau_{k+1} \leftarrow \min\{\tau_k,\max\{\tau_k^{aux},\tau_k/a,\tau_{min}\}$}
    {\If{q>5}{
          $\tau_k^{aux} \leftarrow a\tau_k^{aux}$   
    }
     $\tau_{k+1} \leftarrow \min\{\tau_k^{aux},10\tau_k\}$
    }
    
\end{algorithm}

\subsubsection{Level-bundle methods}

The level bundle method (LBM) is the opposite of PBM and the trust-region methods discussed below \citep{Frangioni2020}. Instead of maximizing the predicted descent of the model $(U^{up}+\phi(Capa^{up})) - (U^{(k+1)}+\hat{\phi}^{(k)}(Capa^{(k+1)}))$ subject to a constraint or a penalty on $\|Capa^{(k+1)} - Capa^{up}\|$, LBM determines a minimum descent target and minimizes the deviation from the stability center needed to reach it (\cite{lemarechal_new_1995}). The stabilized MP for LBM is given by:

\begin{align} \label{eq:lvl}
    \min_{\{Capa,Exp\}} \frac{1}{2}\|Exp-Exp^{up}\|_2^2  \quad s.t. \ \ U(Capa) + \hat{\phi}^{(k)}(Capa) \leq \ell_k, \ \  Capa_{y,i}^{up} = \sum_{y'\in Y_y} Exp_{y',i} \ \forall y,i
\end{align}

The level parameter $\ell_k$ sets the maximum MP objective for the next iterate. It determines level sets of the cutting plane algorithm, i.e., decisions for capacity expansion for which the algorithm estimates function values below $\ell_k$.

The level set is empty, and the problem is infeasible if the objective of the MP is strictly above the prescribed level for all feasible values of $Exp$. This information is still valuable since $\ell_k$ is a lower bound to the overall problem in this case. Therefore, one way to update the lower bound is to set it equal to $\ell_k$ if an iteration resulted in an empty level set \citep{Ackooij2014}. Alternatively, one can update the lower bound by solving the original and stabilized MP in each iteration. One potential advantage of LBM over PBM is that the dynamic parameter $\ell_k$ is in the same order of magnitude as the objective functions, which facilitates setting it to appropriate values.  

As an extension, \cite{Oliveira2016} propose a doubly stabilized bundle method that combines PBM and LBM. This method automatically chooses between proximal and level iterates. The latter occurs when the level set constraint is binding, while proximal iterations find reasonable solutions within a level set instead. The method has the advantage that the strategy to set $\tau_k$ can factor in the dual variable of the level constraint, but due to the method's complexity, we did not include it in our subsequent benchmark.

\paragraph{Parameter strategies}

We devise two different parameter strategies for LBM. The first one strategy, LBM-1, is a straightforward strategy using a fixed weight, $\beta$, to compute a convex combination of the current upper and lower bounds of the problem. We make this adjustment in every iteration, regardless of whether a serious or a null step preceded the adjustment. 
Alg. \ref{alg:lbm1} formalizes this simple mechanism based on \cite{Frangioni2020}. As discussed above, we compute the lower bound by solving the unstabilized MP in every iteration. 

\begin{algorithm} \label{alg:lbm1}
    \small
    \caption{Dynamic adjustment of level parameter $\ell_k$ for LBM-1}
    \DontPrintSemicolon
    \textbf{\textit{Step 6.1 (Initialization)}}:\;
    initialize hyperparameter $\beta$ \;
    \textbf{\textit{Step 6.2 (Define level parameter)}}:\;
    $\ell_{k+1} = \beta\zeta_{low} + (1-\beta)\zeta_{up}$ \;
\end{algorithm}

The second strategy, LBM-2, originates from \cite{brannlund_1995}. 
Two aspects of this strategy are more refined: First, it only sets the level to a convex combination of the upper and lower bounds if the previous iteration performed a serious step ($q>0$). For two, we only increase the level bound in a null step if the level constraint has been binding as indicated by a non-negative dual variable of the constraint $\rho_k$. These refinements ensure that level descent parameter $v^{\ell}_k$ is kept constant for a short sequence of null steps before it is eventually adjusted (cf. \cite{Oliveira2016} for additional reasoning on the approach). The hyperparameters for LBM-2 are $\beta$ and $\mu_{max}$.

\begin{algorithm} \label{alg:lbm2}
    \small
    \caption{Dynamic adjustment of level parameter $\ell_k$ for LBM-2}
    \DontPrintSemicolon
    \textbf{\textit{Step 6.1 (Initialization)}}:\;
    initialize hyperparameters: weight parameter $\beta$ and maximum dual threshold $\mu_{max}$ \;
    \textbf{}
    \textbf{\textit{Step 6.2 (Obtain dual of level set constraint)}}:\;
    Let $\rho_k$ be the current dual solution belonging to constraint $U(Capa) + \hat{\phi}(Capa) \leq \ell_k : \rho_k$\;
    Set $\mu_k = 1+\rho_k$\; 
    \textbf{\textit{Step 6.3 (Adjust level descent)}}: \;
    \eIf{q>0}{$v^{\ell}_{k+1} = \min\{v^{\ell}_{k},(1-\beta)\cdot (\zeta^{up}-\zeta_{low})\}$}{
    \If{$\mu_k > \mu_{max}$}{$v^{\ell}_{k+1} = \beta\cdot v^{\ell}_{k}$}
    }
    \textbf{\textit{Step 6.4 (Set level)}}: \;
    $\ell_{k+1} = \zeta^{up} - v^{\ell}_{k+1}$
\end{algorithm}

\subsubsection{Trust-region bundle methods} \label{qtr}

The final bundle method stabilizes the algorithm using a \textit{quadratic trust-region} (QTR). The next iterate is constrained to be in a Euclidean hypersphere around the current stability center defined by a certain radius. The quadratically constrained stabilized MP is given by:

\begin{align} \label{eq:tr}
    \min_{\{Capa,Exp\}} U(Capa) + \hat{\phi}^{(k)}(Capa) \quad s.t. \ \ \|Exp-Exp^{up}\|_2^2 \leq \gamma_k, \ \ Capa_{y,i}^{up} = \sum_{y'\in Y_y} Exp_{y',i} \ \forall y,i
\end{align}

Note that PBM is the Lagrangian relaxation of the QTR (\cite{HiriartUrruty1993}), highlighting the theoretical equivalence of the methods.

The existing literature does not use QTRs, but commonly applies a linear alternative using the $l_{\infty}$-norm named Boxstep method \citep{Marsten1975}. Although it keeps the stabilized MP linear, \citet{Frangioni2014} show it is inferior to other bundle methods.

\paragraph{Parameter strategies}

The dynamic parameter for the trust-region methods, QTR and Boxstep, is the radius $\gamma_k$ constraining deviations from the current stability center measured with the $l_2$-norm and the $l_{\infty}$-norm, respectively. 

For QTR, we describe the corresponding steps of the strategy in Alg. \ref{alg:qtr}. Rather than directly adjusting the radius itself, we introduce a scaling factor $\psi_k$, that defines the radius $\gamma_k$ relative to the $l_1$-norm of the current stability center. As a result, the radius changes in response either when the dynamic parameter changes or when a serious step changes the stability center.

\begin{algorithm} \label{alg:qtr}
    \small
    \caption{Dynamic adjustment of radius $\gamma_k$}
    \DontPrintSemicolon
    \textbf{\textit{Step 6.1 (Initialization)}}:\;
    initialize hyperparameters $\psi_1,\psi_{min},\kappa\in(0,1],\varepsilon>0$ \;
    \textbf{\textit{Step 6.2 (Reduction test)}}:\;
    \eIf{$\left|1-\frac{U^{(k)}_{unstab}}{U^{(k)}}\right|<\varepsilon$}{
    $\psi_{k+1} \leftarrow \min\{\psi_{min},\psi_k\cdot \kappa\}$        
    }{$\psi_{k+1} \leftarrow \psi_k$}
    \textbf{\textit{Step 6.3 (Update radius)}} \;
    $\gamma_{k+1} \leftarrow (\psi_{k+1})^2 \cdot \|Exp^{up}\|_1^2$
\end{algorithm}

\begin{figure}[!htp]
	\centering
		\includegraphics[scale=0.75]{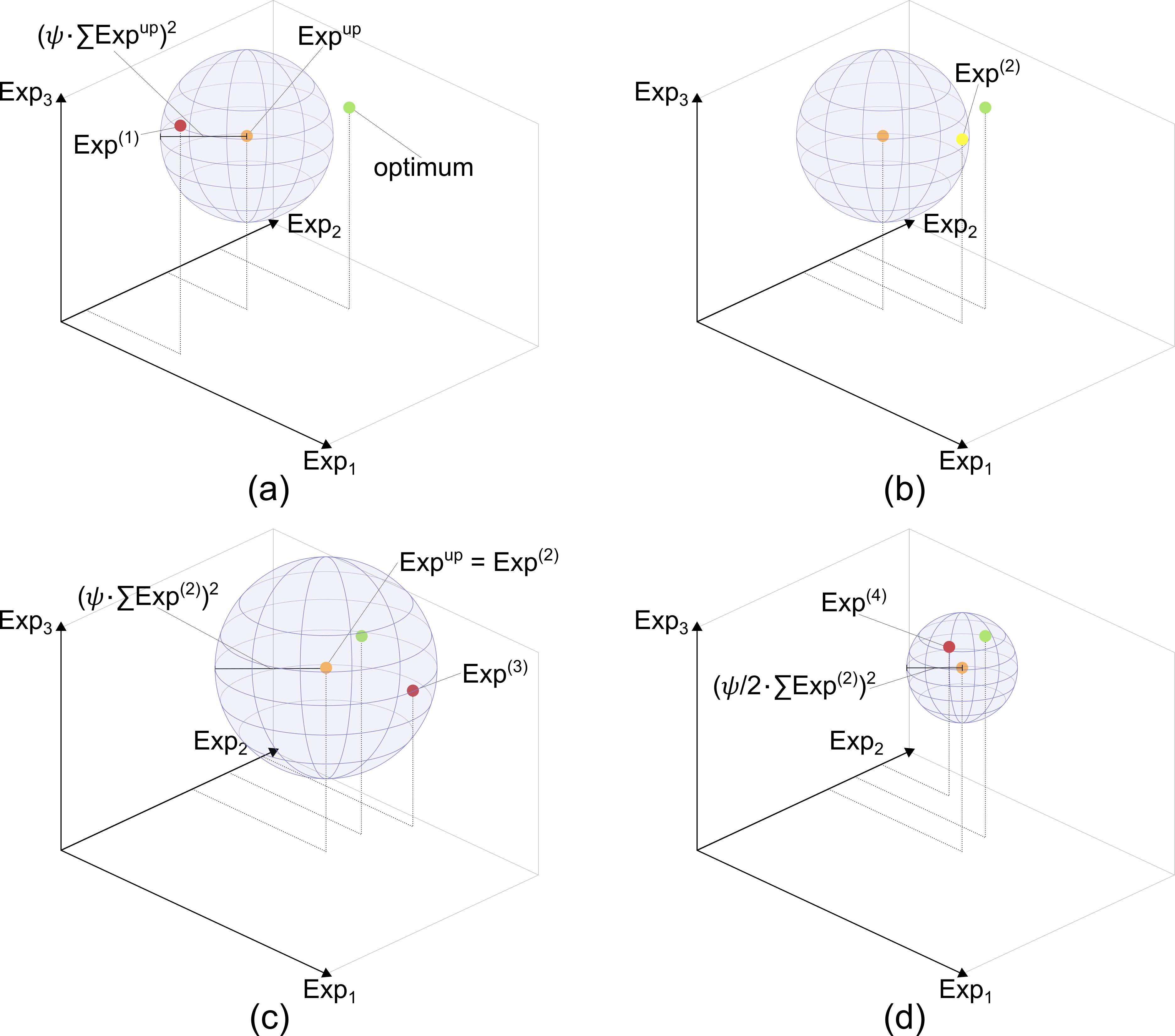}
	\caption{Dynamic trust-region adjustment}
	\label{fig:trustReg}
\end{figure}

We illustrate the stabilization based on four consecutive iterations for a stylized example in Fig. \ref{fig:trustReg}. It includes three unbounded expansion variables, restricted by a spherical trust-region defined by the constraint of Eq. \ref{eq:tr}. 

The current stability center $Exp^{up}$ is indicated in orange, and the solution at the current iteration $Exp^{(k+1)}$ in yellow or red, depending on whether it improves the current best or not. The minimizer of $U(Capa) + \phi(Capa)$ is shown in green. In the first iteration, the distance between the origin and the initial stability center determines the sphere's radius. The first iteration in (a) results in an objective above the current best, but the trust-region is binding, as determined in Step 6.2, and remains unchanged. As the cutting plane model improves locally, the second iteration (b) results in a serious step. Re-centering the trust-region for the third iteration leads to a radius adjustment due to the new $l_1$-norm of the stability center. The third iteration in (c) results in a null step, but the trust-region is no longer binding, and the factor $\psi_k$ is reduced by factor $\kappa$. The fourth iteration with a smaller trust-region in (d) does not improve the solution.

Since a small radius will result in numerical issues, we impose a minimum on the dynamic parameter $\psi_k$, $\psi_{min}$. 
The hyperparameters include the starting radius scaling factor $\psi_1$, the minimum scaling factor $\psi_{min}$, and the reduction factor $\kappa$. In addition, there is a threshold  $\varepsilon$ to check if the trust-region is non-binding since the original and stabilized MP will rarely have the same objective, even if the trust-region is non-binding due to numerical inaccuracies. 

Adding a rule for increasing the radius again is also conceivable because a small starting radius scaling factor $\psi_1$ can result in slow convergence. However, this implies additional hyperparameters, and in practice, it is much easier to correct a small starting factor if the first iterations show slow progress.

For the Boxstep method, we test a range of different box sizes, which we also denote by $\psi$. We define the box size relative to the values for the stability center, so the current values impact the absolute boundaries imposed by the trust-region, resulting in a noncubical box (\cite{Marsten1975}). In order to avoid problems with values of zero, we also impose a minimum upper limit.

\section{Additional refinements} \label{algRef}

In addition to stabilization, there are various methods to improve the poor convergence of the basic BD algorithm. Most refinements, however, aim at problems with a complex combinatorial MP rather than large SPs, like the two-stage planning problem examined in this paper.

One refined already mentioned in Alg. \ref{alg:stdBD} and \ref{alg:bundle},  is deleting non-binding cuts after a predefined number of iterations. The literature on bundle methods refers to such strategies as  $\beta$-strategies \citep{Frangioni2002}. In the algorithms, $\delta_{k,y,s}$ tracks in which iteration a cut was created or binding the last time is initialized based on the number of the current iteration $k$. Step 4 loops over all previous cuts, updating $\delta_{k,y,s}$ and deleting cuts that were not created or binding within the last $\eta$ iterations.

Furthermore, the outlined algorithm is subject to several practical refinements. All capacity variables have an upper limit that is far beyond any reasonable value to pre-limit the solution space of the MP. We scale all equations according to the methodology outlined in \citep{Goeke2020a} to improve numerical properties. If scaling does not suffice, small values of $Capa_{y,i}^{(k)}$ in the SPs are rounded off to zero. The same applies to $\lambda_{y,s,i}^{(k)}$ when adding cuts to the MP. For storage, we enforce an upper and lower bound on the ratio between storage and energy capacity that reflects technical restrictions and confines the solution space.

The following sections discuss additional refinements to the stabilized BD in Alg. \ref{alg:bundle}, namely, inexact oracles, valid inequalities, and parallelization.

\subsection{Inexact oracles} \label{inex}

A large body of literature addresses bundle methods and oracle that return inexact information of the SPs \citep{de_oliveira_convex_2014,Ackooij2014,deOliviera2020,}. Although it is generally possible that oracles cannot obtain exact information, in our application, the oracles deliberately only return inexact information $\{\phi_{y,s}(Capa),\lambda_{y,s,i}\}$  to reduce the computational time of solving the SPs.

\begin{figure}[!htp]
	\centering
		\includegraphics[scale=0.50]{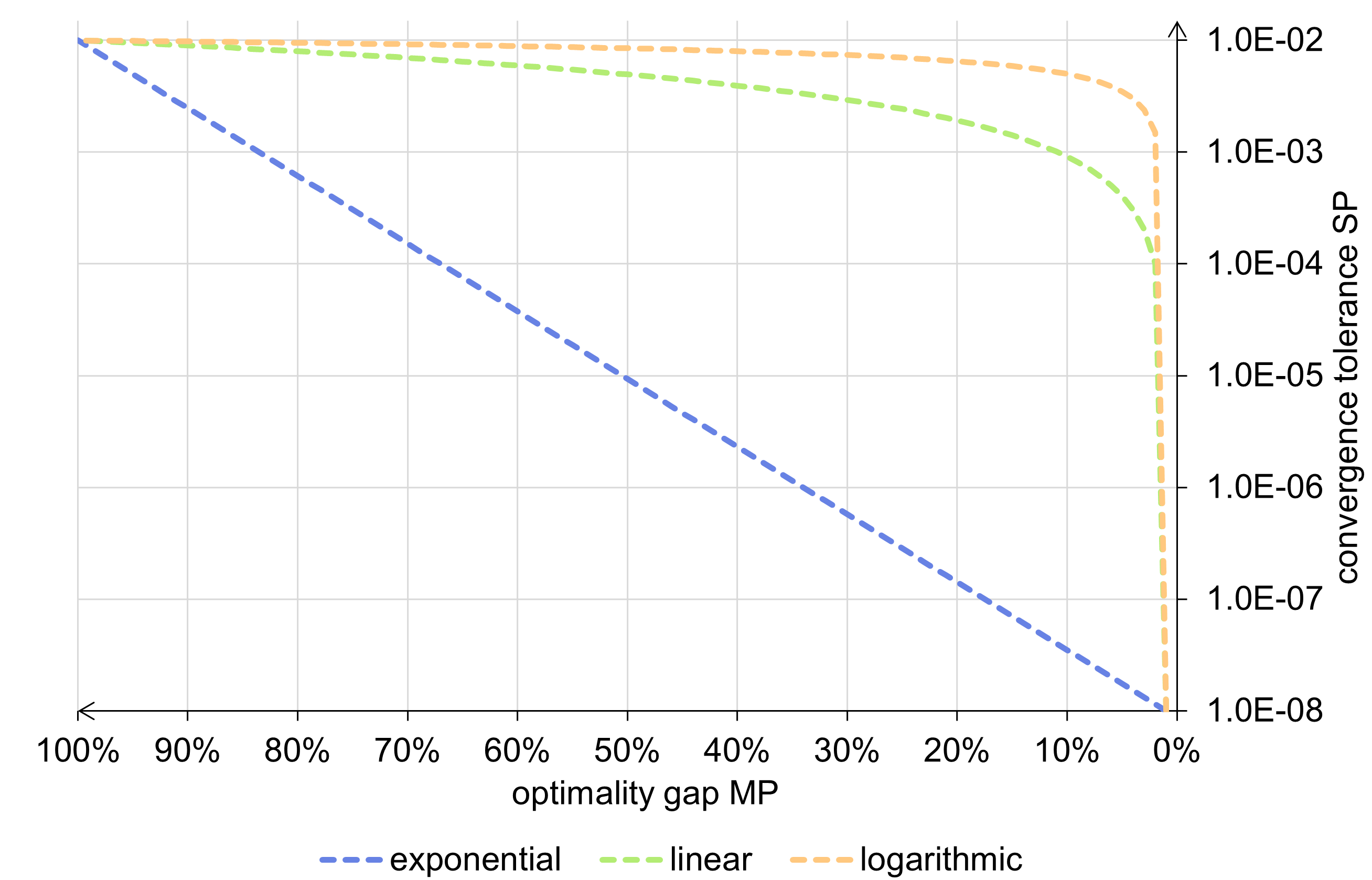}
	\caption{Interpolation methods for convergence tolerance of SPs}
	\label{fig:inexact}
\end{figure}

In practice, we implement this feature by deactivating the crossover step of the Barrier algorithm. As a result, the SPs will return a feasible interior point instead of the optimal basic solution. In addition, we test an asymptotically exact oracle strategy by decreasing the Barrier convergence tolerance as the optimality gap closes, eventually reaching the desired optimality level \citep{Zakeri2000,Ackooij2014}. This strategy faces a trade-off: On the one hand, increasing the tolerance reduces the computation time of the SPs, but on the other hand, it result in less accurate cuts and a potential increase of iterations. Against this background, Fig. \ref{fig:inexact} presents three different interpolation methods for this process. In all cases, the algorithm starts at a tolerance of $1e^{-2}$ and decreases to $1e^{-8}$, but the reduction follows an exponential, linear, or logarithmic curve. Alg. \ref{alg:bundle} can easily implement this approach without further changes.

Beyond the implemented approaches, on-demand oracle accuracy could improve the algorithm further \citep{de_oliveira_level_2014}. With this method, interaction with the solver enables an assessment if the next iterate from Step 2 of Alg. \ref{alg:bundle} will result in a serious step. However, pratical implementations in the context of two-stage stochastic problems focus on cases where the use of multi-cuts significantly impact the MP's computation time. In our problem, this is not the case since the number of scenarios is small compared to the constraints in the MP and even with multi-cuts the computation time of the MP is insignificant \cite{wolf_applying_2014,fabian_computational_2013}.

\subsection{Valid inequalities} \label{vi}

\textit{Valid inequalities} (VI) are constraints added to the MP that are redundant in the closed formulation of the problem but effectively reduce the solution space of the MP \citep{Cordeau2006}. Thus, valid inequalities aim to enhance the convergence of BD at the expense of increasing the MP and its solution time. Often, valid inequalities are derived from the SPs and avoid their infeasibility. 

In the examined two-stage problem, the MP does not include an energy balance and will frequently underinvest in generation capacity, particularly in early iterations. To address this, \citet{Wang2013} add one scenario's operational constraints, including the energy balance, to the MP. In our case, this approach is not viable due to the size of the individual scenarios and corresponding SPs. Instead, we add aggregated operational constraints in Eqs. \ref{eq:5a} and \ref{eq:5b} to the MP for each scenario.
{\small
\begin{subequations}
\begin{alignat}{4} 
& \sum_{t \in T} dem_{y,s,t} &&  \, \leq \, && \sum_{i \in I^{ge}} Gen_{y,s} && \,\,\,\, \forall y \in Y, s \in S \label{eq:5a}  \\
& Gen_{y,s,i} && \, \leq \, && \sum_{t \in T} cf_{s,i,t} \cdot Capa_{y,i}^{ge} &&   \,\,\,\, \forall y \in Y, s \in S, i \in I^{ge} \label{eq:5b} 
\end{alignat}
\end{subequations}
}
The first equation \ref{eq:5a} enforces a yearly energy balance, setting a lower bound for generation that is binding when hourly generation patterns precisely match demand without relying on storage methods that incur energy losses. The second equation \ref{eq:5b} establishes a connection between the lower limit on generation and the capacity variables. As discussed in section \ref{form}, the problem formulation in the paper is still stylized and omits the multiple regions in the case study connected by transmission infrastructure. Therefore, the applied valid inequalities are aggregated by region, assuming lossless and unrestricted transmission as a lower bound.

As highlighted in section \ref{form}, the problem formulation in the paper remains stylized and does not account for the multiple regions connected by transmission infrastructure in the case study. Consequently, the applied valid inequalities also aggregate demand and generation by region, assuming lossless and unrestricted transmission as a lower bound.

\subsection{Parallelization} \label{parallel}

As shown in Alg. \ref{alg:stdBD} and \ref{alg:bundle}, we parallelize solving the SPs. As a result, we can also solve the original MP in every iteration when solving the SPs to obtain a lower bound for the algorithm at no computational costs since the MP solves much faster than any of the SPs. This computation of the lower bounds is common for LBM methods with simple MPs relative to the SPs \citep{Oliveira2016}. It also enables us to implement the parameter strategy for the quadratic trust region.

The efficiency of parallelization is the observed speed-up divided by the number of distributed computation nodes, which corresponds to the number of SPs in our case. In the ideal case, parallelization has no overhead, and all SPs require the same time to converge, resulting in an efficiency of one.

\section{Case study} \label{case}

This section introduces the case-study used to test and benchmark the algorithms introduced previously. The ambition of the case-study is not to achieve highly accurate system planning but to provide a plausible problem covering all critical elements of renewable energy systems. The public repository linked in the Supplementary material section provides All files to re-produce the case-study.

The applied planning model includes two expansion years, 2030 and 2040, focusing on the power sector. The planned system is not subject to any emission constraints in 2030 but must fully decarbonize by 2040 to cover a diverse range of system setups. Spatially, the model covers four distinct regions: France, Belgium, the Netherlands, and Germany. 

\begin{figure}[!htp]
	\centering
		\includegraphics[scale=0.50]{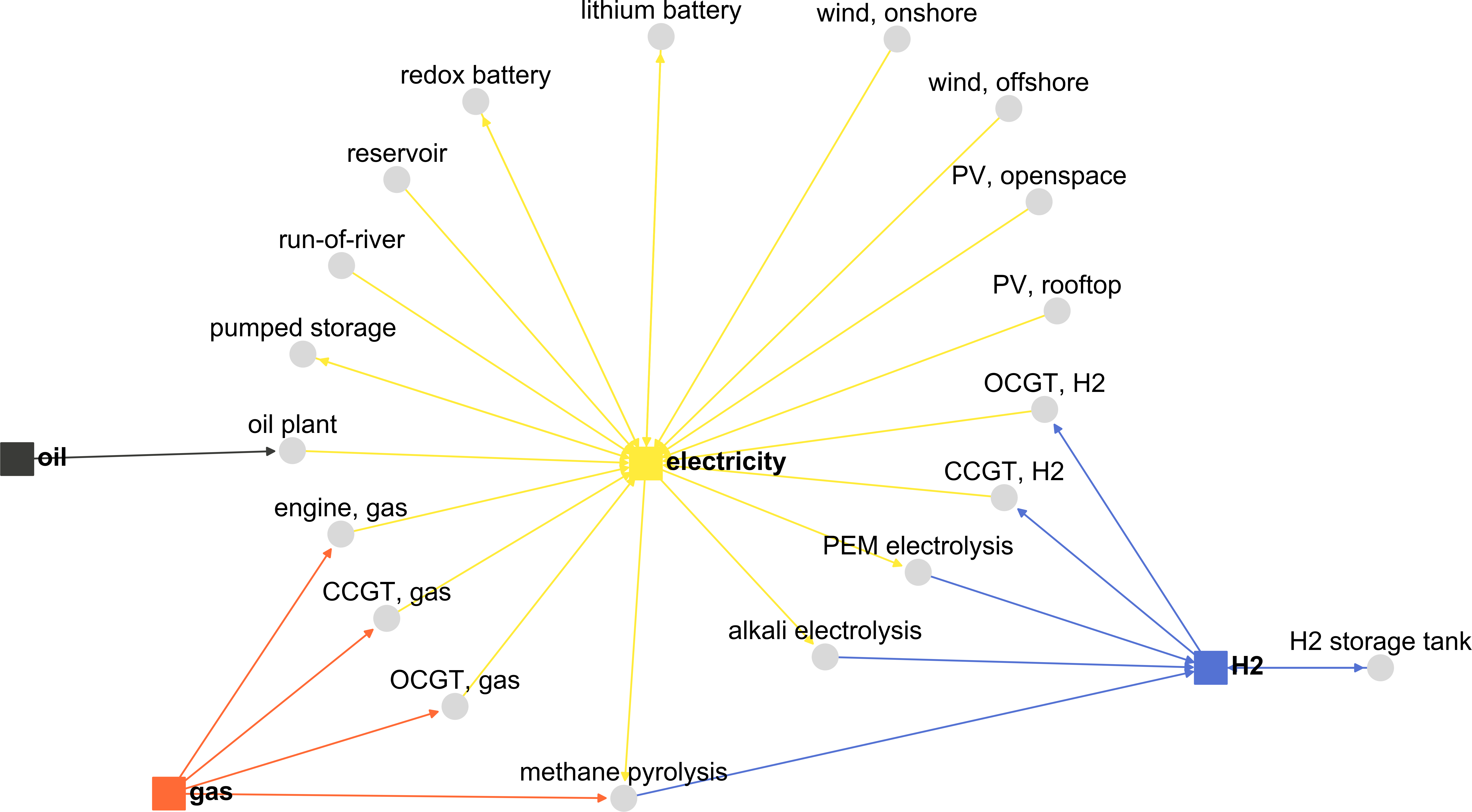}
	\caption{Graph of model elements}
	\label{fig:energyGraph}
\end{figure}

Figure \ref{fig:energyGraph} introduces all considered technologies, depicted as gray circles, and their interaction with energy carriers, depicted as colored squares. Exogenous demand in the model is limited to the carrier electricity modeled at an hourly resolution for the entire year. Hydrogen uses a daily resolution; oil and gas a yearly resolution. 

In the graph, entering edges of technologies refer to their input carriers; outgoing edges relate to outputs. For example, the generation technology electrolysis uses electricity as an input to generate hydrogen. Storage technologies, like pumped storage, have an entering and an outgoing edge to represent charging and discharging. Beyond pumped storage, the model includes batteries for short-term and hydrogen tanks for long-term energy storage. Including long-term storage is critical since it is key for balancing seasonal fluctuations in decarbonized systems. 

Wind and PV capacity are subject to an upper limit that restricts expansion. Pumped storage, run-of-river, and hydro reservoirs have capacities fixed to today's levels. Furthermore, hydro reservoirs are modeled as storage with an exogenous inflow instead of flexible charging. Beyond the listed technologies, the model also decides on the expansion and operation of transmission to exchange electricity and hydrogen between regions.

\section{Results} \label{bench}

The first two sections compare the different stabilization methods and additional refinements. Subsequently, we benchmark selected configurations to an off-the-shelf solver for the deterministic equivalent. This benchmark solves the case study while varying the number of stochastic scenarios from 2 to 16. The MP and each SP run on independent computing cluster nodes, each with four cores and 16 GB of memory. Core models can differ between runs and nodes. Accordingly, the number of nodes amounts to one for the MP plus twice the scenario number for the SPs since the case study covers two expansion years, 2030 and 2040. We re-ran unexpected performance outliers to ensure the robustness of the results. For solving problems, each node deploys the Barrier implementation of Gurobi 10.1 without crossover and the \textit{NumericFocus} parameter set to zero. For BD, we delete unused cuts after 20 iterations and choose a convergence tolerance $\epsilon$ of 0.2\%. The algorithm sometimes fails to converge at smaller tolerances due to rounding. All deployed variations of BD are implemented in the version of the AnyMOD.jl modeling framework linked in the Supplementary material. 

\subsection{Stabilization} \label{res1}

Tab. \ref{tab:tabStab} gives an overview of the tested configurations for the bundle methods and strategies introduced in section \ref{ref2}. For LBM-5, PBM-1, and PBM-2, we copied the configurations proposed in \citet{Oliveira2016}. For all other methods, we chose the configurations ourselves. The initial benchmark covers all permutations of values in the table below. Accordingly, we tested six configurations for LBM-1, 15 for LBM-2, and 27 for PBM-1 and PBM-2 each. For box-step the test covers 4 configurations; for QTR 7. Since this totals 80, initially, each configuration is only tested once using a model with four different scenarios. These tests do not use any refinements other than the respective stabilization.

\begin{table}[]
\begin{tabular}{l|l|c}
Method                  & Parameter strategy & Configurations  \\ \hline
\multirow{2}{*}{Level bundle method} & LBM-1           & $\beta \in \{0.125,0.25,0.375,0.5,0.625,0.75,0.875\}$         \\ \cline{2-3}
                              & LBM-2      & $\beta \in \{0.2,0.5,0.7\}, \mu_{max} \in \{1,5,10,50,100\}$  \\ \hline
\multirow{2}{*}{Proximal bundle method} & PBM-1           & $a \in \{2,4,5\},\tau_1 \in \{1, 5, 10\}$             \\ 
                              & PBM-2          & $\tau_{min} \in \{10^{-6},10^{-5},10^{-5}\}$        \\  \hline  
\multicolumn{2}{l|}{\multirow{2}{*}{Box-step}} & $\psi \in \{2.5\%, 5\%, 7.5\%, 10\%\}$ \\ 
  \multicolumn{2}{l|}{} & upper limit at least 0.5 GW  \\  \hline     
\multicolumn{2}{l|}{\multirow{2}{*}{Quadratic trust-region}} & $\psi_1 \in \{1, 0.5, 0.1, 0.05, 0.01, 0.005, 0.001\}$ \\ 
  \multicolumn{2}{l|}{} & $\psi_{min} = 1 \cdot 10^{-6}, \kappa = 0.5, \varepsilon = 5 \cdot 10^{-4}$ 
\end{tabular}
\caption{Tested configurations for different stabilization methods and Parameter strategies}
\label{tab:tabStab}
\end{table}

\begin{figure}[!htp]
	\centering
		\includegraphics[scale=0.50]{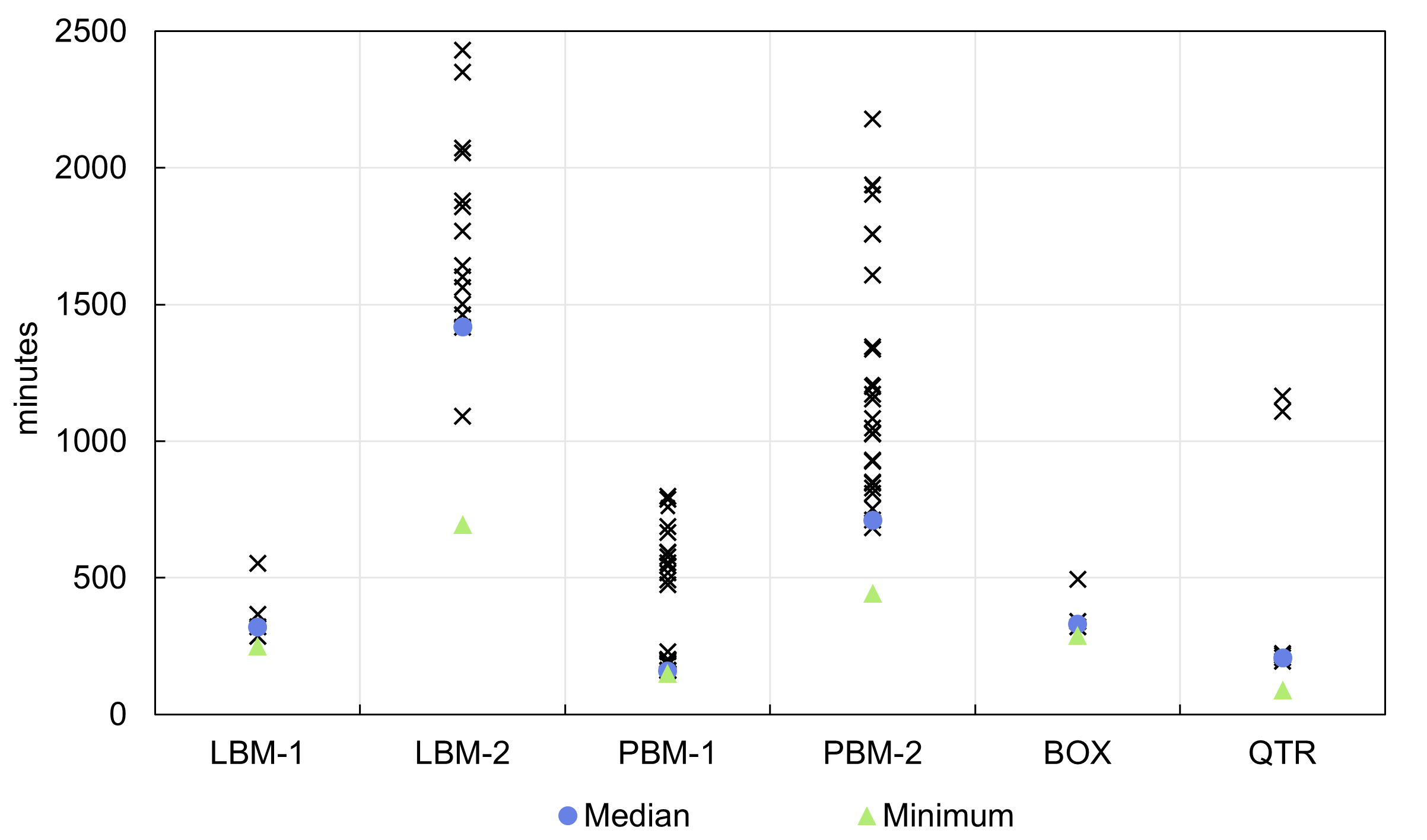}
	\caption{Comparison of configurations for stabilization methods}
	\label{fig:generalMethod}
\end{figure}

Fig. \ref{fig:generalMethod} compares the algorithm's run time when deploying the different configurations. Colored symbols indicate the median and average performance for each group of strategies. Without any stabilization, the computation time amounts to 4'600 minutes, so all methods can substantially accelerate the algorithm, but effects greatly vary by method and configuration. As expected, stabilization improves performance by substantially reducing the number of iterations. For instance, the QTR with $\psi_1 = 0.01$ reduces the number of iterations to 30 from 660 without stabilization. The adverse impact on the MP is negligible since it accounts for less than 0.4\% of the total run-time in both cases. Stabilization even decreases the average time per MP since quick convergence reduces the number of cuts added to the MP that continuously increase its size. 

Overall, LBM-2 and PBM-2 are less performative than the other methods. For these methods and PBM-1, it is also challenging to identify a relation between the configuration and performance. For instance, in the case of PBM-1 and PBM-2, the first and second best configurations of both Parameter strategies do not share any values. In addition, the spread of computation times is significant. For LBM-1, performance generally improves for smaller values of $\beta$ that tighten the stabilization, but the overall spread of results is small; results for BOX are similar. QTR also gives similar results until the initial radius $\psi_1$ reaches a threshold of $0.005$. At this point, the radius is too tight, constraining the MP too much and slowing convergence.

In the following, we reduce the configurations to the two best for each method according to Tab. \ref{tab:selCon}.  
\begin{table}[]
\begin{tabular}{l|l|c|c}
Method & Parameter strategy & \multicolumn{2}{c}{Configurations}  \\ \hline
\multirow{2}{*}{Level bundle method} & LBM-1           & \cellcolor{blue!25} $\lambda = 0.125$ & $\lambda = 0.25$          \\ \cline{2-4}
                              & LBM-1       & $\lambda = 0.5, \mu_{max} = 100$ & $\lambda = 0.7, \mu_{max} = 100$ \\ \hline
\multirow{2}{*}{Proximal bundle method} & PBM-1           & \cellcolor{blue!25} $a =4.0,\tau_1 = 1.0, \tau_{min} = 10^{-6}$ & $a =5.0,\tau_1 = 10.0, \tau_{min} = 10^{-5}$              \\ \cline{2-4}
                              & PBM-2    &       $a =5.0,\tau_1 = 5.0, \tau_{min} = 10^{-5}$ & $a =4.0,\tau_1 = 1.0, \tau_{min} = 10^{-3}$ \\  \hline  
\multicolumn{2}{l|}{Box-step} & $\psi = 7.5\%$ & \cellcolor{blue!25} $\psi = 10\%$ \\ \hline     
\multicolumn{2}{l|}{Quadratic trust-region} & $\psi_1 = 0.05$ &  \cellcolor{blue!25} $\psi_1 = 0.01$ 
\end{tabular}
\caption{Two best configurations for each stabilization method and Parameter strategies}
\label{tab:selCon}
\end{table}

In the previous tests, we solve the model for the most probable scenarios with a daily instead of an hourly resolution to obtain an initial solution for the stabilization. In principle, stabilization benefits from a reasonable initial solution, but there is no strict relationship between the quality of the initial solution and convergence time. Since the algorithm cannot identify an optimal solution in the beginning, poorer starting solutions can even result in a shorter computing time. In addition, obtaining better initial solutions will also increase pre-processing time.

\begin{figure}[!htp]
	\centering
		\includegraphics[scale=0.50]{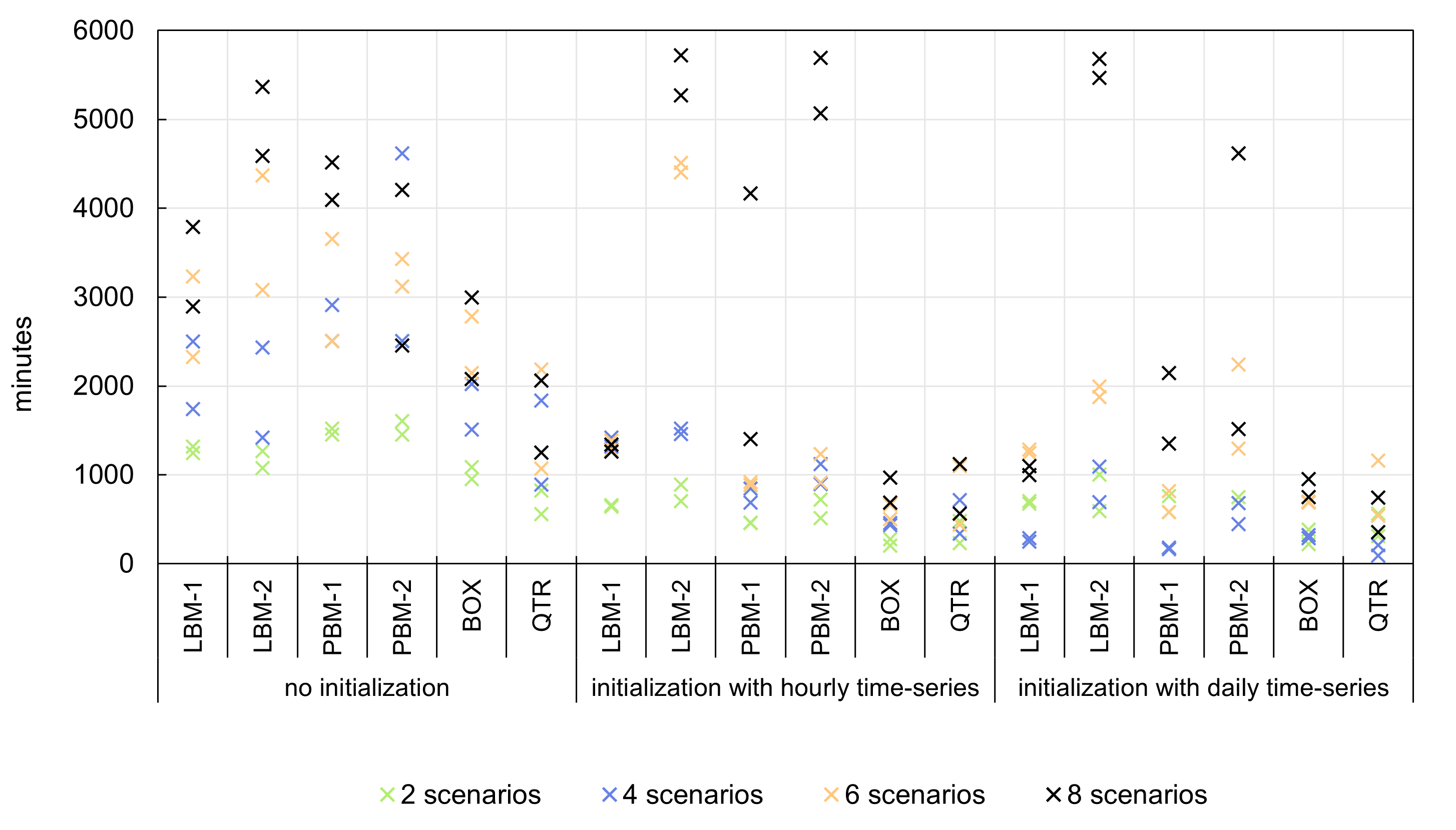}
	\caption{Comparison of initialization strategies}
	\label{fig:iniStab}
\end{figure}

Fig. \ref{fig:iniStab} compares initialization with the most probable scenarios at a daily resolution on the right, to no initialization on the left, and initialization with the most probable scenarios but an hourly resolution in the middle. We run the tests for the selected subset of configurations, varying the number of scenarios from 2 to 8. Again, the results show a large spread; in some instances, initialization significantly impacts a specific configuration's performance. Such interactions are expected since the starting point is an important determinant for the stabilized iteration. However, predicting the effect initialization will have on a specific configuration is impossible. On average, stabilization with the most probable scenario and a daily resolution performs best. Therefore, we use this configuration in the following.  

\subsection{Additional refinements}

In the next step, we analyze the impact other refinements have on the algorithm's performance. 

\begin{figure}[!htp]
	\centering
		\includegraphics[scale=0.50]{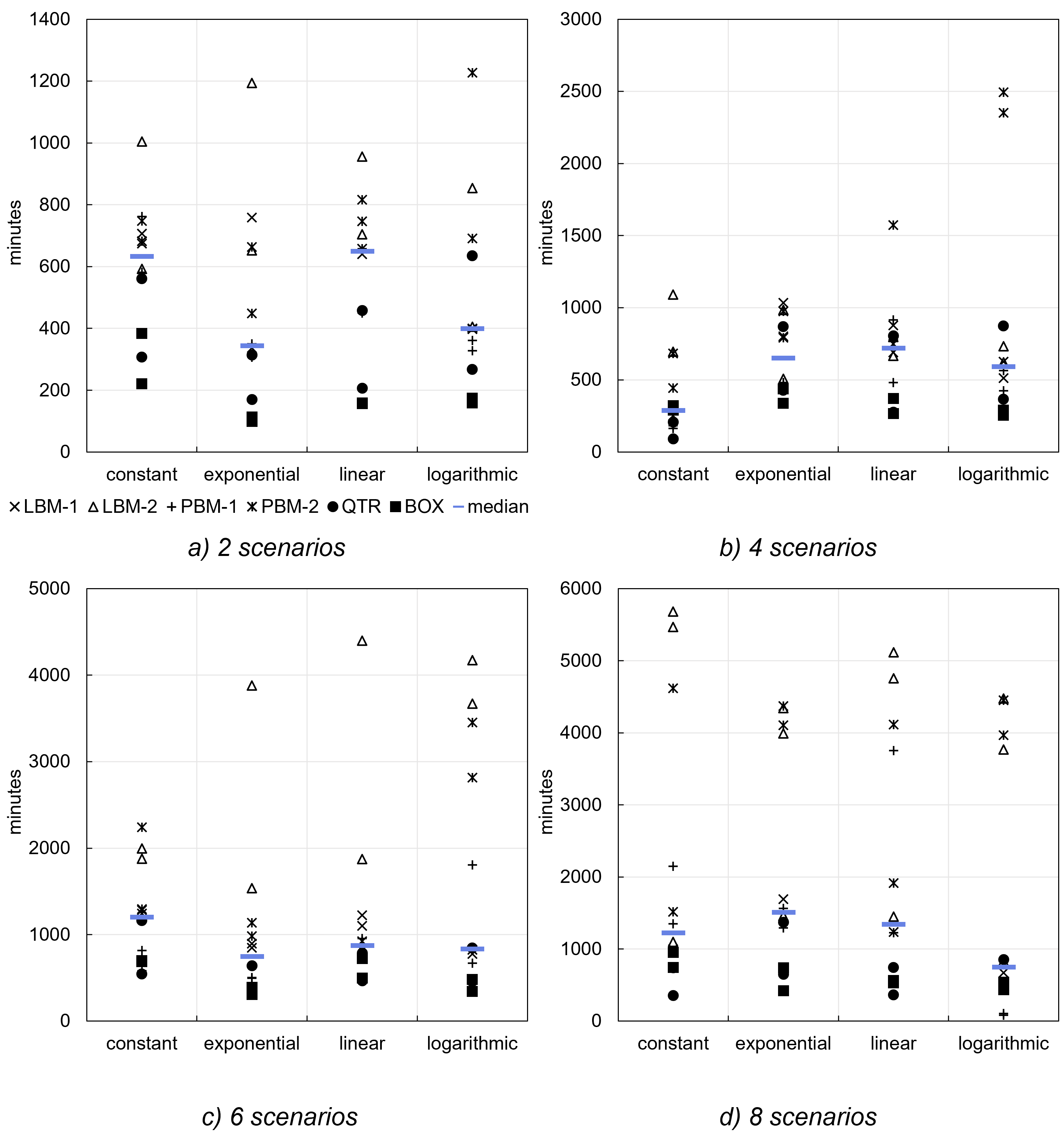}
	\caption{Impact of inexact cuts}
	\label{fig:inexCuts}
\end{figure}

First, Fig. \ref{fig:inexCuts} shows the impact the inexact-cuts strategies introduced in section \ref{inex} have on performance. As described, the four compared cases only vary the strategy for the convergence tolerance of the SPs. Since all skip the crossover phase of the Barrier algorithm that moves from an interior point to a basic solution, they all use inexact cuts. Without this feature, computation times are unreasonably slow, preventing an efficient benchmark.

The results show no clear benefit of inexact-cut strategies on performance. For specific scenarios and configurations, inexact cuts do reduce computation time. For instance, with six and eight scenarios, the logarithmic strategies achieve a median reduction of computation time by 33\% and never significantly increase computation time. However, no clear trend is observed for the other scenario setups and inexact-cut strategies. For some configurations, non-constant strategies can even significantly increase computation time. What is particularly surprising is that in some cases, for example, in all tests with four scenarios, the average time per subproblem increases. In theory, inexact cuts reduce SP time but may increase the number of iterations. Due to the inconclusive results, we presume our benchmark using constant convergence tolerances.

\begin{figure}[!htp]
	\centering
		\includegraphics[scale=0.50]{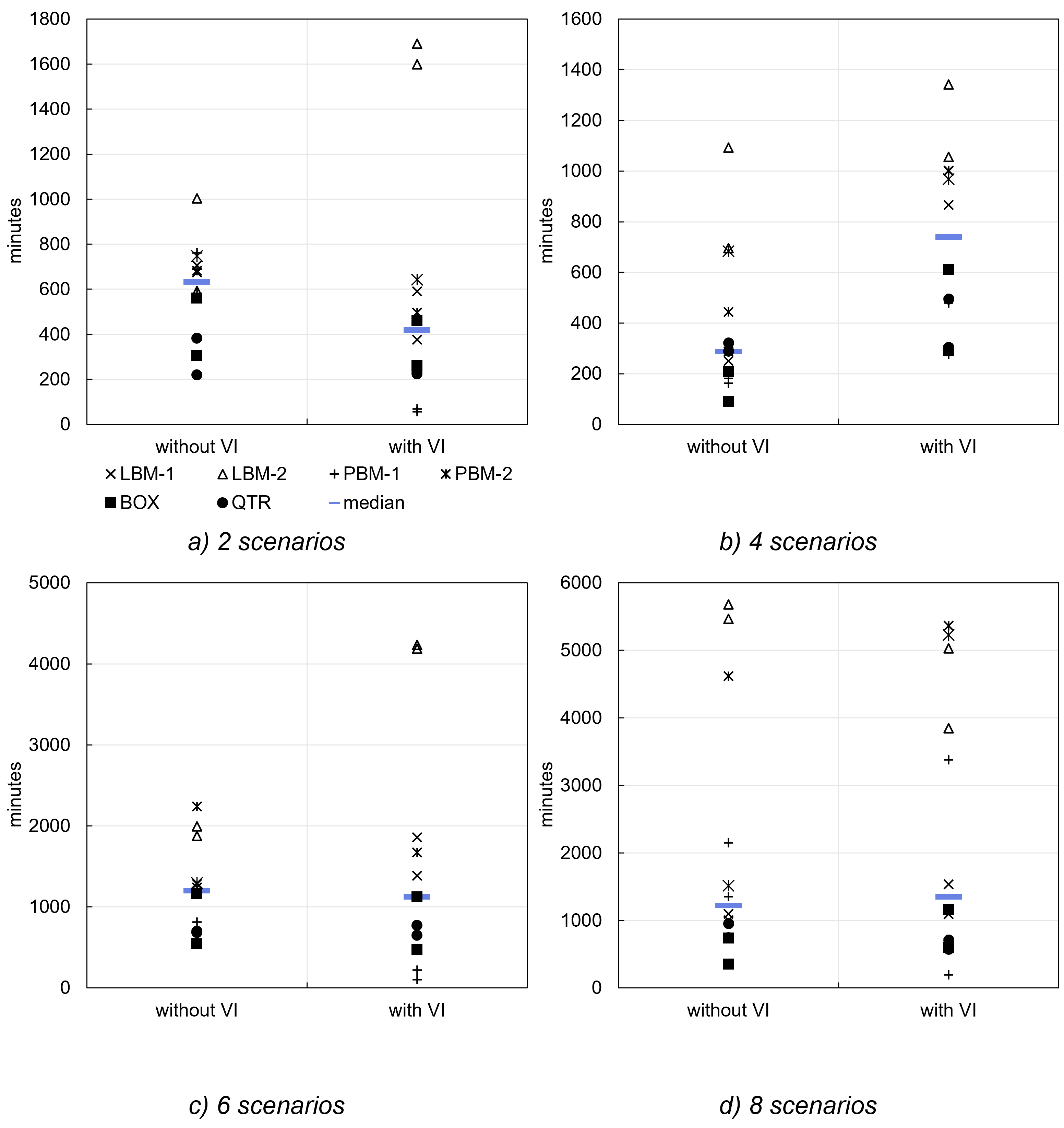}
	\caption{Impact of valid inequalities}
	\label{fig:vi}
\end{figure}

Next, we compare the impact of valid inequalities in the MP as introduced in \ref{vi} on performance. Again, Fig. \ref{fig:vi} compares performance for a different number of scenarios and configurations of stabilization. Overall, valid inequalities increase the spread of solution times without a robust overall reduction. Similar to stabilization, the increase of computation time for the MP is negligible, but we do not observe a consistent reduction of iterations. Even for specific stabilization methods, there is no identifiable trend. For instance, with four scenarios and QTR with a radius $\psi_1 = 0.05$, valid inequalities more than double computation time, but at six scenarios, reduce time by 16\%. Therefore, we will not deploy valid inequalities in the subsequent benchmarks.

\begin{figure}[!htp]
	\centering
		\includegraphics[scale=0.50]{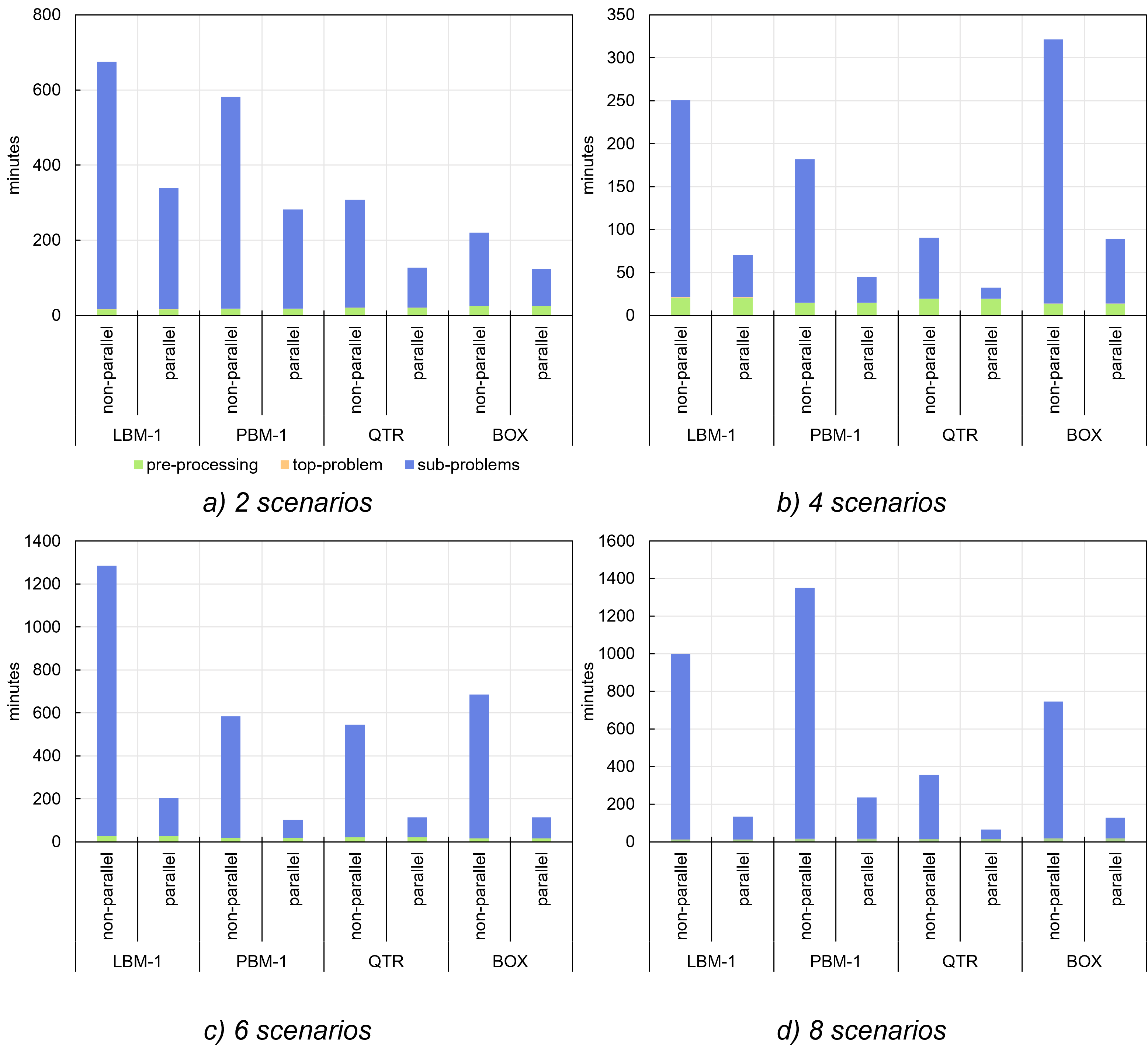}
	\caption{Impact of parallelization}
	\label{fig:para}
\end{figure}

Finally, Fig. \ref{fig:para} shows the benefit of parallelization. As discussed in section \ref{parallel}, parallelization does not strictly change the algorithm itself but distributes the SPs across different nodes to solve them in parallel. If the overhead is negligible, the refinement will strictly decrease computation time but requires additional computational resources.

Across all tests, the efficiency of parallelization varies between 43\% and 68\% with an average of 57\%. The overhead is negligible, and different solve times of the SPs account for almost the entire efficiency loss. To a certain extent, differences in the underlying years can explain the fluctuations of solve times, but mostly, they are random. As the number of scenarios increases, there is no added risk of outliers delaying the algorithm; efficiency stays constant, and the speed-up increases with the number of nodes.

\subsection{Benchmark with the deterministic equivalent}

In the final section, we perform a concluding evaluation of the refined Benders algorithm and compare it to solving the deterministic equivalent with an off-the-shelf solver. Such comparisons are always at risk of cherry-picking and comparing the best algorithm setup determined in an extensive series of tests against a single test with an off-the-shelf solver. 

In the first step, we reduce the variations of stabilization methods and configurations from ten to four, indicated by the blue shade in Tab. \ref{tab:selCon}. This selection builds on the extensive testing with two, four, six, and eight scenarios in the previous sections. Apart from that, all four reference cases use daily initializations, constant convergence tolerances of the SPs, no valid inequalities, and parallelization.

\begin{figure}[!htp]
	\centering
		\includegraphics[scale=0.50]{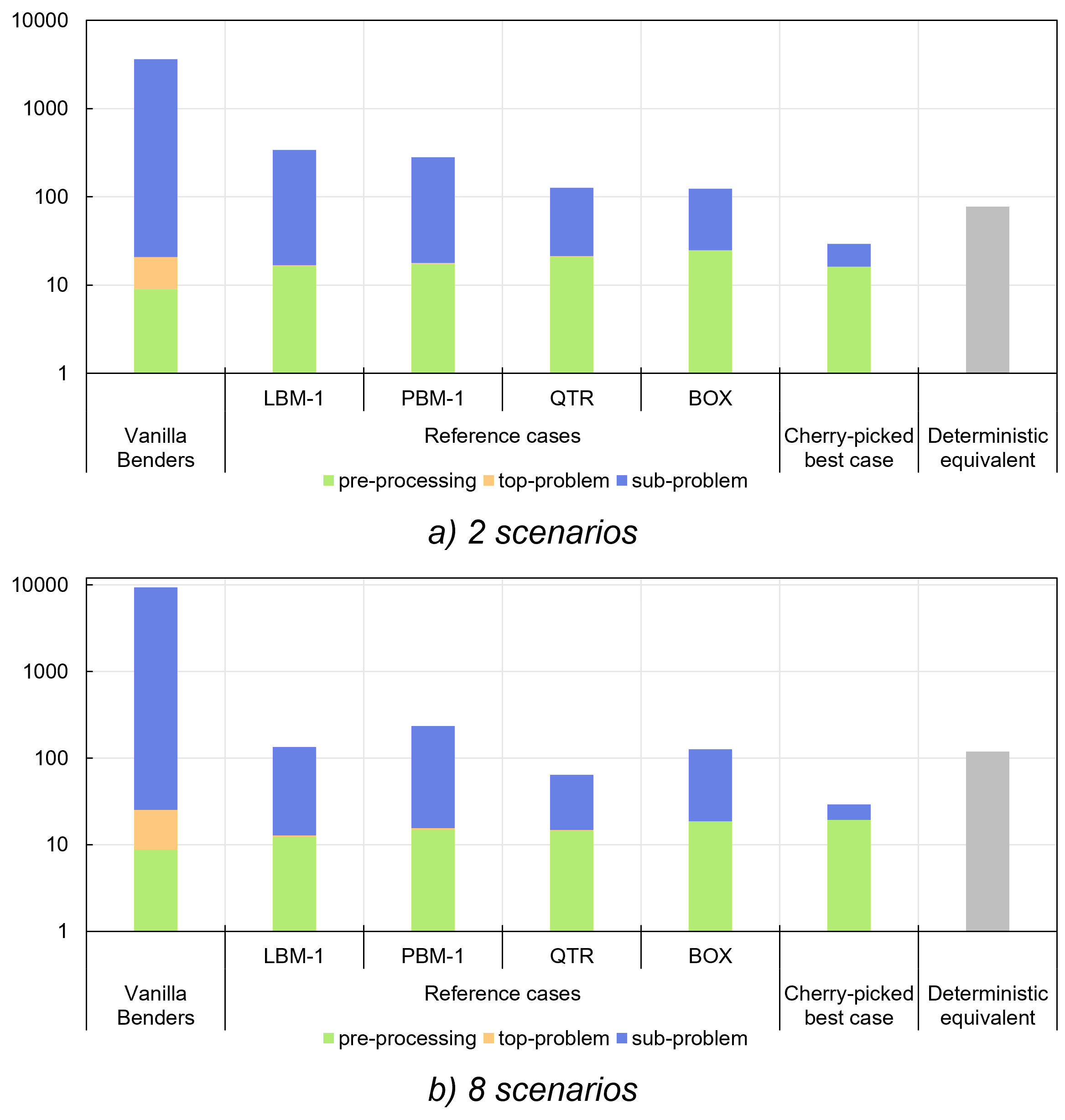}
	\caption{Overall benchmark of refined algorithm}
	\label{fig:bench1}
\end{figure}

Overall, the graph shows a massive speed-up of BD around a factor of 100 thanks to stabilization and parallelization. For two scenarios, the cherry-picked setup takes 29.4 minutes, outperforming the deterministic equivalent with 77.3 minutes, but the more robust reference cases still require between 132.2 (BOX) and 339.0 minutes (LBM-1). With eight scenarios, the reference cases take between 38.2 minutes (PBM-1) and 134.2 minutes (LBM-1) on average, outperforming the deterministic equivalent that requires 119 minutes.

\begin{figure}[!htp]
	\centering
		\includegraphics[scale=0.50]{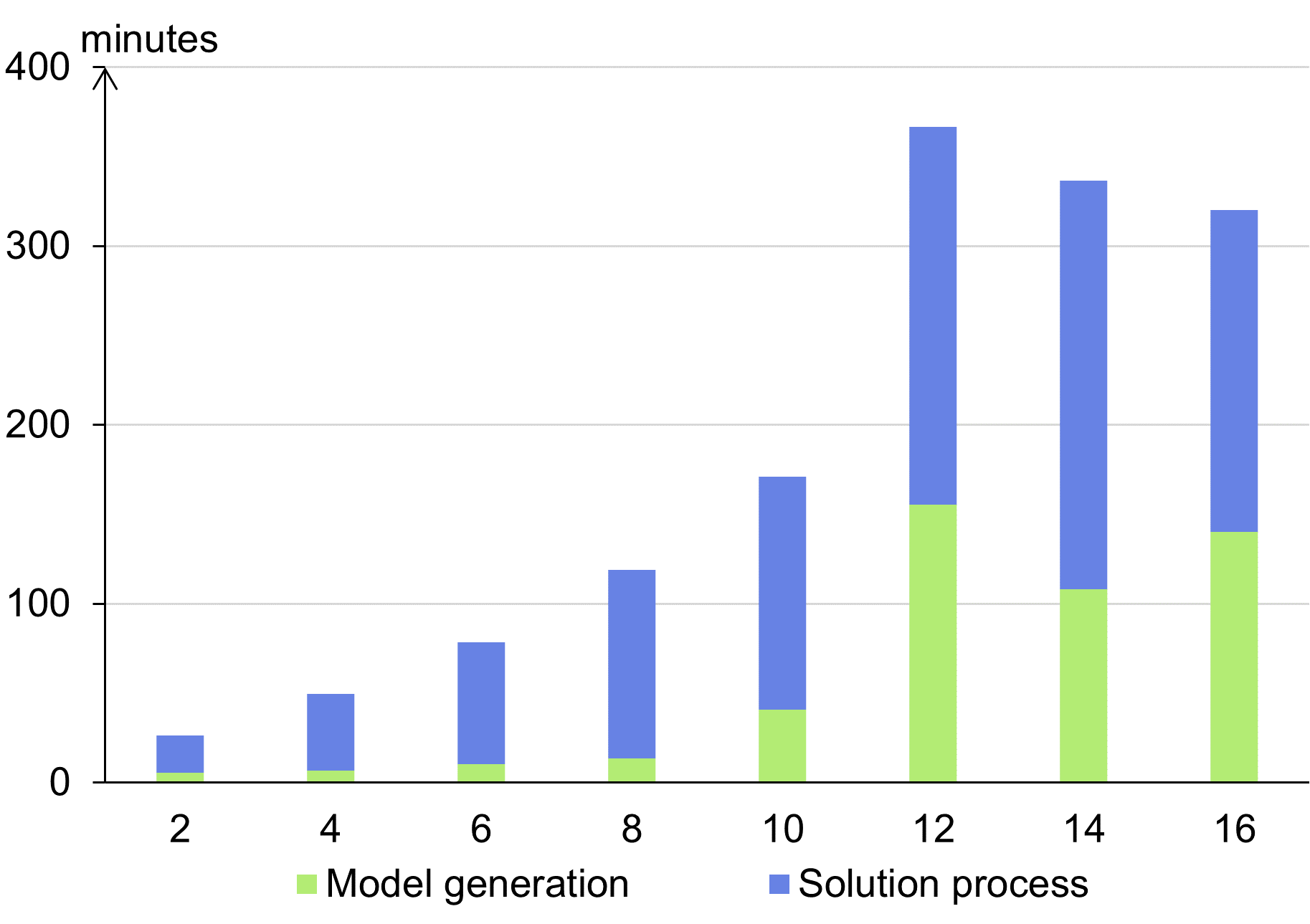}
	\caption{Out-of-sample benchmark of reference cases}
	\label{fig:bench2}
\end{figure}

\begin{figure}[!htp]
	\centering
		\includegraphics[scale=0.50]{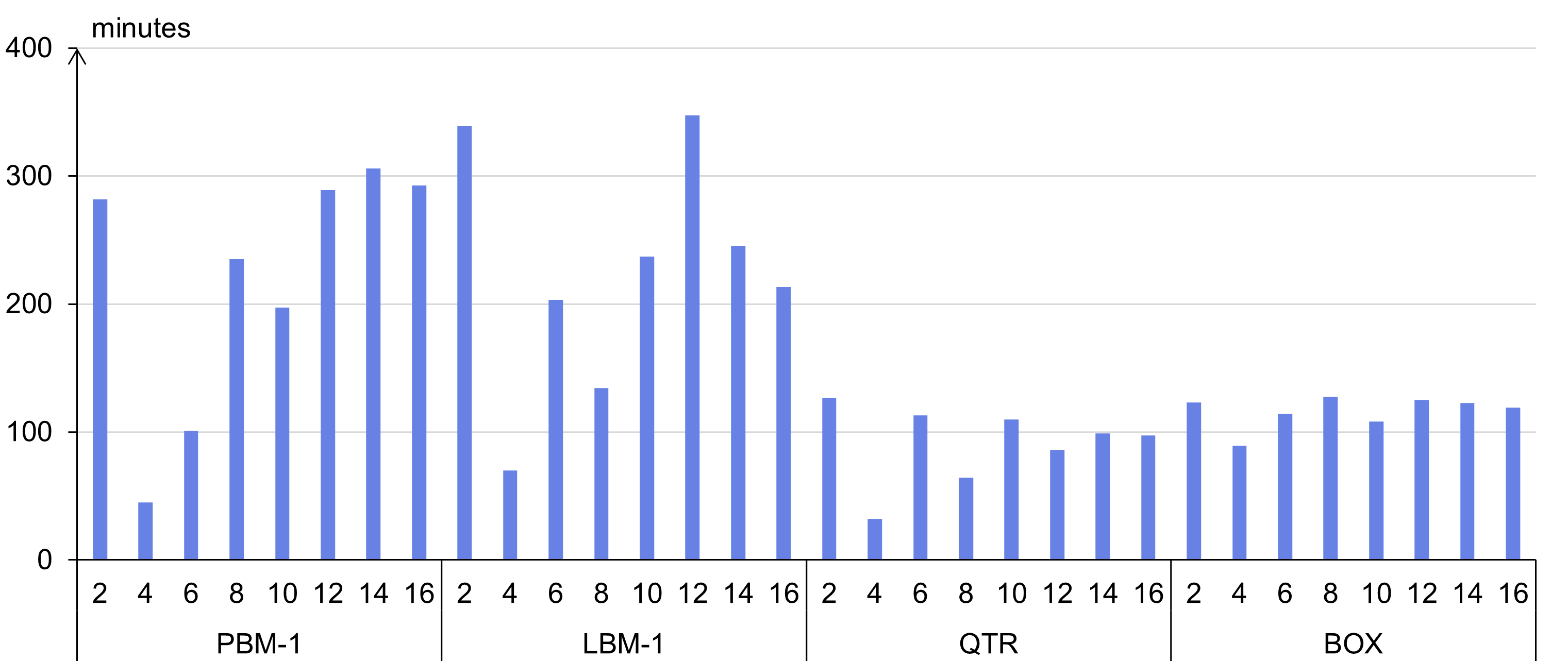}
	\caption{Performance of deterministic equivalent}
	\label{fig:bench3}
\end{figure}

However, the comparison is still biased since we selected the reference cases based on their performance with two to eight scenarios. Therefore, Fig. \ref{fig:bench2} defined reference case "out-of-sample" with up to 16 scenarios and compares their performance to the corresponding deterministic equivalents in Fig. \ref{fig:bench3}. 

Computation time for the deterministic equivalent increases exponentially, although several substantial occur for 14 and 16 scenarios. The graph shows a significant fluctuation of results for the selected proximal and level setups, and methods do not consistently outperform the deterministic equivalent at a certain threshold. For QTR and Boxstep, computation times barely fluctuate, and the selected setups reliably outperform the deterministic equivalent. QTR and Boxstep also showed little fluctuations when comparing different parameter configurations in section \ref{res1}, while the proximal method fluctuated substantially. There is no significant difference in solution quality between the refined BD algorithm and the deterministic equivalent. On average, objective values for the deterministic equivalent are 0.5\permil smaller.

This comparison has limitations since the deterministic equivalent and the parallelized algorithm utilize computational resources differently. In our case, the deterministic equivalent runs on a single node with 24 cores and 16 GB memory each, the most potent node available to us, while BD utilizes up to 33 nodes with 4 cores and 16 GB memory each. Beyond computational resources, performance can vary for other model setups beyond the range tested in this paper. Since the benchmark is not universally valid, the focus should not be on solution times but on scaling behavior.

\section{Conclusion}

This paper applied BD to two-stage stochastic problems for energy planning under climate uncertainty, a critical problem for designing renewable energy systems. To improve performance, we modify the standard algorithm and adapt various refinements to the problem's characteristics---a simple continuous MP, and few but large SPs. These refinements include inexact cuts, valid inequalities, and parallelization, but the focus is on stabilization methods to mitigate the oscillating behavior of BD. We test a broad range of bundle methods and configurations, including a quadratic trust-region, which is theoretically equivalent to other quadratic methods but novel for continuous problems. 

In a quantitative case-study, all stabilization methods can significantly reduce computation time compared to the standard algorithm. However, the quadratic trust-region and the non-quadratic box-step method have practical advantages. Firstly, these methods have proven to be robust because their performance is less sensitive to the numeric configuration. As a result, they require less fine-tuning to a specific problem to achieve good results. Second, the practical implementation of these methods is more straightforward.

The developed algorithm can also benefit from inexact cuts and valid inequalities, but there are also adverse effectand no clear trend. On the other hand, parallelization consistently improves performance by a factor corresponding to the number of SPs divided by two, which translates to a parallelization efficiency of 50\%. Overall, the refined algorithm accelerates the vanilla BD algorithm by a factor of 100. Thanks to parallelization, its computation time remains constant when the number of scenarios increases. Cherry-picked setups can already outperform off-the-shelves solvers for the deterministic equivalent for two scenarios; more general algorithm setups start to achieve this around six scenarios, but note that performance can vary for other model setups beyond the range tested in this paper.

There are further conceivable improvements to the algorithm described in this paper. To list them, we again follow the categorization of enhancements introduced in
\citet{Rahmaniani2017}:
\begin{itemize}
  \item \textbf{Solution procedure}: Considering most of the run-time is spent on the SPs, solving approximations instead of the original SPs can further increase performance. Not solving SPs to optimality in the current algorithm already exploits that valid cuts do not require exact solutions. Conceivable approximation methods include surrogate models based on machine learning, geometric interpolation, or using a reduced temporal resolution \citep{Neumann2021,Mazzi2021,Goeke2021a}. However, approximations face two challenges: First, the high dimensionality of inputs to the SP corresponds to the number of complicating capacity variables. Second, using asymptotically inexact oracles also requires a revision of the stabilization method. The use of on-demand accuracy oracles appears promising. Introducing a solver component that allows to skip iterations and replace by them with approximations if they do no meet a predefined descent target is a natural extension of the current algorithm \citep{Ackooij2014}.
  
  \item \textbf{Solution generation}: Although heuristic methods are insufficient for adequate planning, they can greatly narrow the solution space of the MP to improve convergence. This approach also offers synergies with the selection of representative climatic years for the planning problem in the first place. For instance, solving the deterministic problem for specific climatic years, as performed in section \ref{clus}, is not only useful for clustering climatic years, but can also derive lower and upper bounds for capacity variables or system costs.
  
  \item \textbf{Cut generation}: The disproportional impact storage technologies have on the run-time that was discussed at the end of section \ref{res1} suggests the algorithm would benefit from advanced cut generation. However, the size of the SPs still poses an obstacle to this approach. Therefore, refined approaches that only solve a modified but not the original SP are most promising \citep{Papadakos2008,Sherali2013}.

\end{itemize}

Apart from further enhancements, the outlined algorithm enables robust planning of renewable macro-energy systems wit a large number of climatic years, which was the original motivation for this work. How many climatic years to include and how to select these years for adequate system planning is an open question for future research. In addition, future research should extend the optimization under uncertainty to the operational stage in the SPs. Currently, uncertainty is limited to expansion in the MP, and operation in each scenario assumes perfect foresight, but stochastic dual dynamic programming could capture uncertainty at the operational stage as well \citep{Papavasiliou2018}.

Beyond from climate uncertainty, the algorithm has the potential to make a broad range of other analyses in energy planning computationally tractable. In contrast to the deterministic equivalent, the algorithm scales well when adding complexity to the expansion problem since its computation time is negligible in the decomposed algorithm. For instance, the algorithm is suited for planning with endogenous learning, which adds mixed-integer variables to the expansion problem for a piecewise linear approximation of cost curves \citep{Felling2022,Zeyen2022}. Similarly, the algorithm makes a stochastic representation of investment costs or discrete investment decisions in highly resolved models viable. Furthermore, it can efficiently compute near-optimal solutions due to its iterative nature and since oracles remain valid when changing the objective function.

\section{Acknowledgments}

The research leading to these results has received funding from the German Federal Ministry for Economic Affairs and Energy via the project "MODEZEEN" (grant number FKZ 03EI1019D). A special thanks goes to all Julia developers. In addition, we would like to express our gratitude to the reviewers for their insightful feedback on earlier drafts of this paper.

\section*{Supplementary material}

The data and execution scripts for the case-study are available here: \url{https://github.com/leonardgoeke/EuSysMod/releases/tag/bendersPaper}.

\printcredits

\bibliographystyle{apacite}
\bibliography{cas-refs}

\section{Appendix}

\subsection{Sampling of climatic years} \label{clus}

To compare the introduced variations of BD, we solve the model as a two-stage planning problem with multiple scenarios. Each scenario represents a climatic year with different patterns and total levels of electricity demand and capacity factors for wind and PV. For this purpose, we build on corresponding data for 39 climatic years from 1980 to 2018 published in \citet{bloomfield2020}, based on re-analysis data from the NASA's MERRA-2 dataset.

Stochastic programming often reduces $n$ scenarios of historical data, in our case 39 climatic years, to a smaller number $m<n$ that captures the empirical distribution but improves computational tractability \citep{dupa2003,rujeerapaiboon2022}. Scenario reduction methods determine a set of representative scenarios $S$ and a set of corresponding weights $\{p_s:s\in S, \sum_{s\in S} p_s = 1, p_s\geq 0 \}$, such that the weighted scenario average of second stage costs approximates the sample average over all $n$ scenarios.  There is a variety of reduction methods for this purpose ranging from na\"ive sampling to moment matching or clustering methods \citep{kaut2021}. The choice of scenario reduction method and the optimal number of scenarios $m$ for adequate planning of renewable systems is an open question that is beyond the scope of this paper. For a meaningful comparison of solution methods, as conducted in Section \ref{bench} below, it is sufficient to ensure the $m$ selected scenarios are not extremely similar and bias the algorithmic performance as a result.

In this paper, we used a problem-dependent k-mediod algorithm, as proposed by \citet{Mundru2022}, for selecting $m$ representative scenarios and corresponding probability weights. The method is advantageous for large-scale planning problems since it reduces the high-dimensional clustering problem including hourly data on demand and capacity factors across multiple countries to a much smaller clustering problem in $\mathbb{R}_+$, the space of system costs. To perform the clustering, we calculate a matrix of symmetric distances between scenarios, or climatic years, where the $s$-th row and $s^*$-th column is computed according to Eq. \ref{eq:dis}:

\begin{align} 
    d_{s,s^*} = \frac{1}{2}\left[Z_{s,s^*} - Z_{s,s}+Z_{s^*,s}-Z_{s^*,s^*}\right] \label{eq:dis}
\end{align}

where $Z_{s,s^*}$ refers to total system costs when the deterministic problem for scenario $s$ is solved with capacities fixed to results of the deterministic problem for scenario $s^*$. For instance, $Z_{1980,2008}$ corresponds to total system costs, including infeasibility costs, when solving the deterministic model for the climatic year 1980 with capacities originally computed for 2008. The clustering algorithm selects $m$ representative scenarios and assigns each of the $n-m$ scenarios remaining in the sample to one of the representative scenarios, such that the sum of distances between representative scenarios and assigned scenarios as defined above, is minimized.

\begin{figure}[!htp]
    \centering
    \includegraphics[scale=0.1]{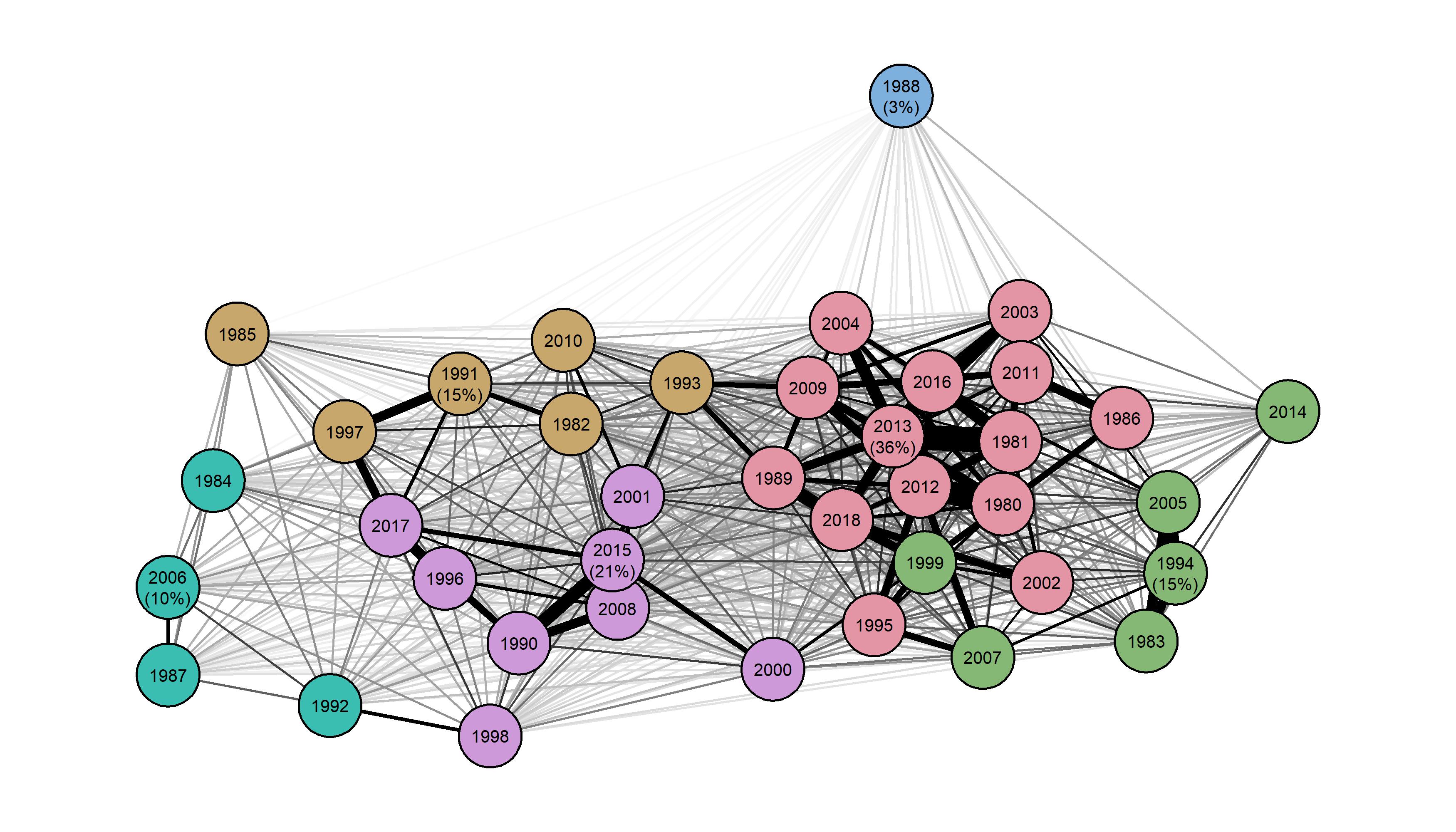}
    \caption{Similarity graph for problem-dependent scenario reduction with $m=6$ clusters}
    \label{fig:pdsr}
\end{figure}

Fig. \ref{fig:pdsr} illustrates the results of the clustering algorithm for $m=6$. Each color represents a group and the node of the representative scenario for this group is equipped with its weight $p_s$, which reflects the relative group size. The color and thickness of connecting lines between the nodes indicate \textit{similarity}, the inverse of the distance $d_{s,s^*}$ as defined above.

\end{document}